\documentclass{article}

\usepackage{shapepar,bm}
\usepackage[dvips,arrow,matrix,ps,color,line]{xy}
\usepackage{theorem,amsfonts,amssymb,amscd}
\usepackage{geometry}
\usepackage{makeidx}

\usepackage{graphicx}
\usepackage{graphics}

\usepackage{epstopdf}

\usepackage{amssymb, graphics,hyperref}

\hypersetup{
    colorlinks = true,
    linkcolor = blue,
    anchorcolor = red,
   citecolor = blue,
    filecolor = red,
    pagecolor = red,
    urlcolor = blue}

\newcommand{\mathsym}[1]{{}}

\usepackage[usenames]{color}
\definecolor{MyLightMagenta}{cmyk}{0.1,0.8,0,0.1}
\definecolor{MyDarkBlue}{rgb}{0.1,0,0.3}
\makeindex

\def\NN{\mathbb N}

\def\bfu{{\mathbf u}}

\def\mod{{\mathrm{mod}}}

\def\ZZ{\mathbb Z}

\def\CC{\mathbb C}

\def\QQ{\mathbb Q}

\def\cocoa{{\hbox{\rm C\kern-.13em o\kern-.07em C\kern-.13em o\kern-.15em A}}}

\def\Dcal{\mathcal D}

\def\Span{{\rm Span}}
\def\hom{{\mathrm{Hom}}}

\def\Ical{\mathcal I}

\def\bfu{{\bf u}}

\def\blamb{{\bm \lambda}}

\def\Pcal{{\mathcal P}}

\def\w2M{\bigwedge^2M}

\def\wM{\bigwedge M}

\def\w{\wedge }
\def\bw{\bigwedge }

\protect
\protect
\protect
\protect

\protect
\protect

\def\sra{\rightarrow}
\def\lra{\longrightarrow}

\def\proof{\noindent{\bf Proof.}\,\,}
\def\qed{{\hfill\vrule height4pt width4pt depth0pt}\medskip}
\def\be{\begin{equation}}
\def\ee{\end{equation}}
\def\bclm{\begin{claim}}
\def\eclm{\end{claim}}
\def\beqn{\begin{eqnarray}}
\def\eeqn{\end{eqnarray}}
\def\beqn*{\begin{eqnarray*}}
\def\eeqn*{\end{eqnarray*}}

\theoremstyle{change}
\theorembodyfont{\rmfamily}

\newtheorem{claim}{}[section]
\global
\global

\def\no@breaks#1{{\def\\{ \ignorespaces}#1}}    % disallow explicit line breaks

  \makeatletter
\def\cleardoublepage{\clearpage\if@twoside \ifodd\c@page\else
\hbox{} \thispagestyle{empty}
\newpage
\if@twocolumn\hbox{}\newpage\fi\fi\fi} \makeatother

\usepackage{eso-pic,graphicx}
\makeatletter
\newcommand\BackgroundPicture[2]{%
  \setlength{\unitlength}{1pt}%
  default \put(0,\strip@pt\paperheight){%
  \parbox[t][\paperheight]{\paperwidth}{%
    \vfill
     \centering \includegraphics[angle=#2, width=15cm, height=15cm,  bb=0 0 150 150]{#1}
    \vfill
}}} %
\makeatother
\usepackage{palatino}

\title{Vertex operators arising from  Linear ODEs\footnotetext{2010 {\sl Mathematics Subject Classification}: 14M15, 15A75, 17B69, 34A35, 34A99.}}

\author{\sc{Letterio Gatto, Parham Salehyan}}

\date{}

\begin{document}

\maketitle

\abstract{\noindent  \small{The Heisenberg Oscillator Algebra  admits irreducible representations both on the
ring $B$ of polynomials in infinitely many indeterminates (the
{\em bosonic representation}) and on a graded-by-{\em charge}  vector
space, the {\em
semi-infinite} exterior power of an infinite-dimensional
$\QQ$-vector space $V$ (the {\em fermionic representation}).
Our main observation is that  $V$ can be realized
as the $\QQ$-vector space generated by the solutions to a generic
linear ODE of {\em infinite order}. Within this framework, the well known {\em boson-fermion} correspondence for the zero charge fermionic space  is a  consequence of the
formula expressing each solution to a linear ODE as a
linear combination of the elements of the universal basis of
solutions. In this paper we extend the picture for linear ODEs of finite order. Vertex operators are defined and fully described in this case.}}

\smallskip
\noindent
{\bf Keywords and Phrases.} {\sl Generic Linear ODEs, Boson-Fermion Correspondence, Vertex Operators}. 

\tableofcontents

\section{Introduction}

\claim{\bf Statement of the Results.} This paper is about {\em vertex operators} related to solutions to generic linear ODE of order $r\in \NN\cup\{\infty\}$. If $r<\infty$,  the main characters of the play are: 

\begin{enumerate}
\item a ring $B_r:=\QQ[e_1,\ldots, e_r]$ of polynomials in the indeterminates 
$(e_1,\ldots, e_r)$ with rational coefficients, the $r$-th {\em bosonic Fock space}; 

\item the polynomial $E_r(t) =1-e_1t+\cdots+(-1)^re_rt^r\in B_r[t]$ 
and the sequence $H_r:=(h_j)_{j\in\ZZ}$  defined by the equality  $
E_r(t)\sum_{n\in\ZZ}h_nt^n=1
$,
holding in the ring of formal power series $B_r[[t]]$; 

\item a vector space $V_r:=\Span_\QQ(u_i)_{i\in\ZZ}$, where  $u_i=\displaystyle{\sum_{n\geq 0} h_{n+i}{t^n\over n!}}$;

\item The $r$-th {\em fermionic Fock space}   of total charge $i\in \ZZ$ defined by
$
F^r_{i}:=\Span_\QQ ( \Phi^r_{i,\blamb})_{ \blamb\in \Pcal_r}
$,
where $\Pcal_r$ is the set of all partitions $\blamb:=(\lambda_1\geq\lambda_2\geq\ldots\geq\lambda_r\geq 0)$ of lenght at most $r$ and 
\be
\Phi_{i,\blamb}^r:=u_{i+\lambda_1}\w u_{i-1+\lambda_2}\w\cdots\w u_{i-r+1+\lambda_r}\w u_{i-r}\w u_{i-r-1}\w \cdots\label{eq:charge_i}
\ee

\item The kernel  $K_r\subseteq B_r[[t]]$  of the generic linear Ordinary Differential Operator 
\be
D^r-e_1D^{r-1}+\cdots+(-1)^re_r:B_r[[t]]\sra B_r[[t]],\label{eq:GLODO}
\ee
where $D$ is the usual formal derivative on formal power series.
\end{enumerate}
It is easily checked (Cf.~\cite{GS1}) that 
$
(u_0, u_{-1},\ldots, u_{1-r})
$
is a $B_r$-basis of $K_r$, i.e. $u_{-i}\in K_r$ for $0\leq i\leq r-1$, and  if $\phi\in K_r$ there are unique $U_i(\phi)\in B_r$ such that  $\phi=\sum_{i=0}^{r-1} U_i(\phi)u_{-i}$ (Cf. \ref{ubas}). This suffices  to prove:

\smallskip
\noindent
{\bf Proposition~\ref{pregmb} (Boson-Fermion correspondence).}  {\em The fermionic Fock space $F^r_0$ is a free $B_r$-module of rank $1$ generated by the {\em vacuum} vector $\Phi^r_0:=u_0\w u_{-1}\w u_{-2}\w\cdots$. More precisely $\Phi^r_{0,\blamb}=\Delta_\blamb(H_r)\Phi^r_0$, where $\Delta_\blamb(H_r):=
\det(h_{\lambda_j-j+i})_{1\leq i,j\leq r}$ is the Schur determinant associated to the partition $\blamb$ and the sequence $H_r$ .}

\smallskip
\noindent
Each $\alpha\in F^r_0$ is a finite $\QQ$-linear combination of typical semi-infinite exterior monomials $\Phi^r_{0,\blamb}$ (Cf. (\ref{eq:charge_i})): the {\em boson-fermion} correspondence $F^r_0\sra B_r$  maps
 $\alpha\mapsto \alpha/\Phi^r_0$, with obvious meaning of the notation (Section~\ref{FFS}). 
 
Since  $B_r=\bigoplus_{\blamb\in\Pcal_r}\QQ\cdot \Delta_\blamb(H_r)$,   the canonical $B_r$-module structure of $F^r_0$  prescribed by Proposition~\ref{pregmb} allows to define  {\em vertex operators} $\Gamma_r(z),\Gamma^\vee_r(z):B_r\sra B_r[[z^{-1},z]]$ through the equalities:
\begin{eqnarray*}
\Gamma_r(z)\Delta_\blamb(H_r)&=&{(\sum_{i\in\ZZ}z^iu_i \w \Phi^r_{-1,\blamb})\otimes_{\QQ}{1_{B_r}} \over \Phi^r_0},\\
\Gamma^\vee_r(z)\Delta_\blamb(H_r)&=&{(\sum_{i\in\ZZ}z^{-i}u^\vee_i \lrcorner \Phi^r_{1,\blamb})\otimes_\QQ {1_{B_r}}\over \Phi^r_0}.
\end{eqnarray*}
Consider the $\QQ$-homomorphisms $G_r(z), G^\vee_r(z):B_r\sra B_r[[z^{-1},z]]$ given by:
\begin{center}
$G_r(z)\Delta_\blamb(H_r):=E_r(z)(\Gamma_r(z)\Delta_\blamb(H_r))\qquad\mathrm{and}\qquad \displaystyle{G^\vee_r(z)\Delta_\blamb(H_r):={z\Gamma_r^\vee(z)\Delta_\blamb(H_r) \over E_r(z)}}.$
\end{center}
One has ({\bf Corollaries~\ref{corhnz} and~\ref{scndlm}}):
 \[
G_r(z)h_n=h_n-{h_{n-1}\over z}\qquad\mathrm{and}\qquad G^\vee_r(z)h_n=\sum_{i\geq 0}{h_{n-i}\over z^i},
\]
so that, in fact, $G_r(z)h_n$ and $G^\vee_r(z)h_n$ belong to $B_r[z^{-1}]$. By extending by $\QQ[z^{-1}]$-linearity, an easy check proves that $G^\vee_r(z)G_r(z)h_n=G_r^\vee(z)G_r(z)h_n=h_n$.  Indeed, it turns out that $G_r(z)\Delta_\blamb(H_r)$ and $G^\vee_r(z)\Delta_\blamb(H_r)$ both belong to $B_r[z^{-1}]$ and that $G_r(z)$ and $G^\vee_r(z)$ are one inverse of the other once they are extended by $\QQ[z^{-1}]$-linearity.  This  is a consequence of the main result of this paper, which is the explicit description of  $\Gamma_r(z)$ and $\Gamma_r^\vee(z)$.

\smallskip
\noindent
{\bf Theorems~\ref{mnth1} and~\ref{mnthm2}.} {\em The operators $G_r(z),G^\vee_r(z)$ commute with taking the Schur determinants, i.e.
\[
G_r(z)\Delta_\blamb(H_r)=\Delta_\blamb(G_r(z)H_r)\qquad \mathrm{and}\qquad G^\vee_r(z)\Delta_\blamb(H_r)=\Delta_\blamb(G^\vee_r(z)H_r),
\]
where $G_r(z)H_r:=(G_r(z)h_n)_{n\in\ZZ}$\,\,\, and \,\,\, $G_r^\vee(z)H_r:=(G^\vee_r(z)h_n)_{n\in\ZZ}$.
Therefore:
\begin{eqnarray}
\Gamma_r(z)\Delta_\blamb(H_r)&=&{1\over E_r(z)}\Delta_\blamb(G_r(z)H_r),\label{eq:int01}\\ \nonumber\\
 \Gamma^\vee_r(z)\Delta_\blamb(H_r)&=&{ E_r(z)\over z}\Delta_\blamb(G^\vee_r(z)H_r).\label{eq:int02}
\end{eqnarray}
}

\smallskip
\noindent
The proof of formula~(\ref{eq:int01}) is based on the following expression holding in the fermionic picture  ({\bf Proposition~\ref{cortj}}):
\be
(\sum_{i\in\ZZ}z^iu_i\w \Phi_{-1,\blamb}^r)\otimes_\QQ {1_{B_r}}={1\over E_r(z)}\exp \left({t\over z}\right)\w \Phi_{-1,\blamb}^r,\label{eq:int03}
\ee

\noindent
while the proof of~(\ref{eq:int02}) is based on the equality

\smallskip
\noindent
%{\bf  Lemma \ref{72}.}  
\be
{(z\sum_{i\in \ZZ} u^\vee_iz^{-i}\lrcorner \Phi^r_{1,\blamb})\otimes_\QQ {1_{B_r}}\over \Phi^r_0}=\left|\matrix{z^{-\lambda_1}&z^{1-\lambda_2}&\cdots&z^{r-1-\lambda_r}\cr h_{\lambda_1+1}&h_{\lambda_2}&\cdots &h_{\lambda_r-r+2}\cr
\vdots&\vdots&\ddots&\vdots\cr
h_{\lambda_1+r-1}&h_{\lambda_2+r-2}&\cdots&h_{\lambda_r}}\right|+\Delta_{(\lambda_1+1,\ldots,\lambda_r+1)}(H_r).\label{eq:int04}
\ee
stated and proven in {\bf Lemma~\ref{72}}.

\medskip
\noindent
Let us stress that the first summand of~(\ref{eq:int04}) is obtained from $\Delta_\blamb(H_r)$ by substituting each entry  $h_{\lambda_i+1-i}$ of its first row    by the monomial $z^{i-1-\lambda_i}$ ($1\leq i\leq r$). 

Formulas~(\ref{eq:int01})--(\ref{eq:int02}) and (\ref{eq:int03})--(\ref{eq:int04}) are new in shape and perspective as they involve in an essential way the use of solutions to a linear ODE of order $r$. They generalize the classical framework of the boson-fermion correspondence, as explained e.g. in~\cite{K, KR},  that arises in the representation theory of the Heisenberg algebra. This last picture is in fact obtained by letting $r$ going to $\infty$.   Our Section~\ref{infty} supplies  a new transparent proof of  the classical expression of the vertex operators  as displayed, e.g.,  in~\cite[p.~54]{KR} or~\cite[Section 4]{Arb}. The proof is based on the fact that  $G_\infty(z), G^\vee_\infty(z)$ are ring homomorphisms (contrarily to  $G_r(z),G^\vee_r(z)$, which are not  for finite $r$) and then can be written as exponentials of a first order differential operator.
It is also probably worth  to mention that  for  $r=\infty$, each $u_{-i}$ ($i\in \ZZ$) remarkably satisfies (in many different ways) the Kadomtsev-Petviashvilii (KP)  equation in the Hirota bilinear form (Cf. Remark~\ref{lodeinfty}).
\claim{\bf Relationships with Schubert Calculus.} Mainly inspired
  by mathematical physics (see e.g.~\cite{Arb,K, KR, KL, neretin}),  the  {\em boson-fermion correspondence} has a number of nice geometrical consequences.  In relatively recent times  its connection with the geometry of the Hilbert schemes of points in the affine plane has been investigated by Nakajima~\cite{nakajima}. More classically,  a system of PDEs encoded in the tensor product of the  {\em vertex operators} $\Gamma_\infty(z),\Gamma_\infty^\vee(z)$ allows to  embed an (infinite dimensional) Universal Grassmann Manifold (\cite[p.~76]{KR}) into the ring  of polynomials in infinitely many indeterminates with complex coefficients. For different applications of vertex operators see also~\cite{PraLa}.

The guiding idea of our investigation was  to look at the boson-fermion correspondence  as an infinite dimensional manifestation of the classical Schubert Calculus for the Grassmann varieties $G:=G(k,\CC^n)$ of $k$-planes in $\CC^n$, as formulated in~\cite{G1} through a derivation on the exterior algebra $\wM$ of a free module of rank $n$. Within this framework,  Schubert cycles in the Chow ring of  $G$ can be seen as a kind of {\em bosonic} counterpart of  generalized wronskians associated to a linear ODE having as coefficients  the Chern classes of the tautological bundle over $G$ (see~\cite{GS1}).

%Schubert Calculus can be played  as in \cite{G1} for all Grassmann varieties $G(k, \CC^n)$ at once,  $1\leq k\leq n$,   through a   derivation on a Grassmann algebra of a free $\ZZ$-module,  indeed a prototype  of a vertex operator.  In that context, Giambelli's formula, a consequence of a suitable {\em integration by parts}, turns out   to be a finite dimensional version of the boson-fermion correspondence.

The explicit bridge between generic linear ODEs  and  Schubert Calculus was established in~\cite{GS1} (see also~\cite{GS2}), where a  Giambelli (or Jacobi-Trudy) formula for {\em generalized wronskians}, associated to a fundamental system of solutions, is proven. It  generalizes the proof of the classical theorem, due to Abel and Liouville,  according which   if the wronskian determinant associated to a fundamental system of solutions to a linear ODE does not vanish at a point,  then it vanishes nowhere. The purely algebraic  treatment  of~\cite{GS1} suggests in a irresistible way to define  generic linear ODEs of infinite (countable) order, with indeterminate coefficients  $(e_1,e_2,\ldots)$ (see also~\cite{GattoAlg}). The basis elements of the fermionic space $F^\infty_0$ can be seen as generalized wronskians associated to a basis of solutions,  which is what  allows to rephrase the classical framework of the boson-fermion correspondence.

Equipping the exterior algebra  of a free module $M$ with a derivation, as in~\cite{G1}, is the same as   endowing $\wM$ with a sequence of endomorphisms satisfying certain Leibniz like rules with respect to the wedge product. Their restriction to  a fixed $k$-th exterior power  has  been studied purely in terms of symmetric functions  by   Laksov and Thorup in~\cite{LakTh1,LakTh2}. They show that the formalism of Schubert Calculus for Grassmann varieties is entirely encoded by the canonical symmetric structure of the $k$th exterior power of a polynomial ring. It  equips the latter with a structure of a  free module of rank $1$ over the ring of symmetric functions in $k$ indeterminates.  Remarkably, at the end of the introduction of~\cite{LakTh1}, the authors claim that  ``In the work of E. Date, M. Jimbo, M. Kashiwara and T. Miwa \cite{DJKM}, Schur functions appear in connection with exterior products ({\em but}) in another context. See also the work of V. G. Kac and A. K. Raina \cite{KR}.''  The present paper  shows, on the contrary,  that  the context is pretty much the same and is based on the classical elegant   interplay (described e.g. in~\cite{MacDonald}) between the complete symmetric functions, the elementary symmetric functions and those which are sum of  powers.

\claim{} The paper is organized as follows.  Section~\ref{preli} collects a few preliminaries, Section~\ref{BFS} defines the $r$-th bosonic Fock space. It will be interpreted as the $\QQ$-polynomial ring generated by the indeterminate coefficients of a generic linear ODE  (Section~\ref{GODE}).

 Fermionic-Fock spaces generated by formal power series related with linear ODEs are introduced in Section~\ref{FFS}. The boson-fermion correspondence for $F^r_0$ is proven, through the application of the {\em universal Cauchy formula} (\ref{eq:univcauchy}). Still in this section the vertex operators $\Gamma_r(z),\Gamma_r^\vee(z)$ are defined. They are discussed in detail in Sections~\ref{vertex} and~\ref{secgammavee}. The final Section~\ref{infty} (re)-computes the shape of the vertex operators associated to linear ODEs of infinite order, re-obtaining classical formulas.

 \medskip
 \noindent
 {\bf Acknowledgment.} {Work sponsored by  FAPESP-Brazil  Processo n. 2012/02869-1,  IBILCE-UNESP, Campus of S\~ao Jos\'e do Rio Preto, and partially by  INDAM-GNSAGA, PRIN ``Geometria sulle Variet\`a Algebriche'' and Politecnico di Torino}. We are grateful to Maxim Kazarian and Inna Scherbak for enligthening discussions,   Simon G. Chiossi and Marcos Jardim for constant warm encouragement.
 
\section{Preliminaries and Notation}\label{preli}
\claim{}\label{convpart} Throughout the paper $\NN$ will denote the non negative integers and $\NN^*$ the positive integers. We denote by $\Pcal$ the set of all partitions, i.e the set of all monotonic non increasing sequences of non negative integers  such that all but finitely many terms are zero. The {\em length} $\ell(\blamb)$ of the partition is the number of its  non zero terms, called {\em parts}.  An arbitrary partition $\blamb$ of length at most $r$ will be written  as $(\lambda_1\geq \lambda_{2}\geq\ldots\geq \lambda_{r})$.  Let  $\Pcal_r(\subseteq \Pcal)$ be the set of all partitions of length at most $r$ and  $a:=(a_i)_{i\in\ZZ}$  any bilateral sequence of elements of some ring $A$. Then  $\Delta_\blamb(a):=\det(a_{\lambda_{j}-j+i})_{1\leq i,j\leq r}\in A$  is the {\em Schur determinant} associated to $\blamb$ and to $a$.

\claim{} Let $A$ be any $\QQ$-algebra and let $A[t]\subseteq A[[t]]$ be the standard inclusion of the polynomial ring into the formal power series. A monic polynomial  $P\in A[t]$  of degree $r$ will be written as:
\be
P:=t^{r}-e_1(P)t^{r-1}+\ldots+(-1)^re_r(P),\qquad e_i(P)\in A.\label{eq:monpol}
\ee
where $(-1)^je_j(P)$ denotes the coefficient of $t^j$.
If $P$ splits in $A$ as a product of $r$ distinct linear factors,   $e_j(P)$  is precisely the $j$-th elementary symmetric polynomial function  in the roots of $P$. 
As $A$ contains the rational numbers, each $\phi\in A[[t]]$ can be written in the form
\be
\phi=\sum_{n\geq 0}a_n{t^n\over n!}, \qquad (a_n\in A).\label{eq:for0}
\ee
If $a_n=a^n$, for some $a\in A$, then $\exp(at)$ stands for the associated {\em exponential} formal power series  $
\exp(at)=\sum_{n\geq 0}a^n{t^n/ n!}.
$
%The {\em formal Laplace transform} $L:A[[t]]\sra A[[t]]$ is defined by
%\[
%L(\phi)=L\left(\sum_{n\geq 0}a_n{t^n\over n!}\right)=\sum_{n\geq 0}a_nt^n.
%\]

\claim{}\label{22} The map $D^j:A[[t]]\sra A[[t]]$, $j\in\NN$, defined by:
\[
D^j\sum_{n\geq n}a_n{t^n\over n!}=\sum_{n\geq 0}a_{n+j}{t^n\over n!},
\]
is  the identity of $A[[t]]$ for $j=0$ while for $j>0$ is the $j$th iterated of $D:=D^1$,  the formal derivative on $A[[t]]$, the unique (infinite) linear extension of the map $D t^n=nt^{n-1}$.
The commutative $A$-subalgebra $A[D]$ of $End_A(A[[t]])$ of  linear ordinary differential operators (ODO) with $A$-coefficients is gotten by  evaluating  each $P\in A[t]$ (monic or not) at $D$. In particular    
\[
\ker P(D):=\{y\in A[[t]]\,|\, D^ry-e_1(P)D^{r-1}y+\ldots+(-1)^re_r(P)y=0\},
\] 
 is the $A$-submodule of $A[[t]]$ of the solutions of the   linear homogeneous ODE $P(D)y=0$.

\section{Bosonic Fock Spaces}\label{BFS}
\claim{} \label{secer} Let $r\in \NN\cup\{\infty\}$ and    $E_r:=(e_i)_{i\in [1,r]\cap \NN}$ be a sequence of $r$ indeterminates over $\QQ$. For instance $E_1=(e_1)$, $E_2=(e_1,e_2)$, \ldots, $E_\infty:=(e_1,e_2,\ldots)$. We write
\be
E_r(t)=1+\sum_{i\in [1,r]\cap \NN}(-1)^ie_it^i\in B_r[[t]]\label{eq:er(t)}
\ee
for  the formal power series having as coefficients the terms of the sequence $E_r$. If $r<\infty$ then $E_r(t)=1-e_1t+e_2t^2+\ldots+(-1)^re_rt^r$.
 Imitating~\cite{KR}, we call $r$th {\em bosonic Fock space} the polynomial $\QQ$-algebra  $B_r:=\QQ[E_r]$. For $s\geq r$ there are obvious $\QQ$-algebra epimorphisms  $B_s\sra B_r$ mapping $e_j\mapsto e_j$  if $j\leq r$ and $e_j$ to  $0$ if  $j>r$.

\claim{}\label{hr}  Let $H_r=(h_i)_{i\in\ZZ}$ and $X_r:=(x_i)_{i\in\NN^*}$ be the sequences of  coefficients of the  formal power series  $H_r(t):=\sum_{n\in \ZZ}h_nt^n$ and $X_r(t)=\sum_{n\geq 1}x_nt^n$ defined through the equalities in $B_r[[t]]$:
\be
H_r(t)= {1\over E_r(t)}=\exp(X_r(t)).\label{eq:for1}
\ee
where $E_r$ is defined by~(\ref{eq:er(t)}). 
The first equality~(\ref{eq:for1}) gives:
\[
\sum_{n\in\ZZ} h_nt^n={1\over 1+\sum_{1\leq i\leq r}(-1)^ie_it^i}=1+ \sum_{n\geq 1}\left(\sum_{1\leq i\leq r}(-1)^ie_it^i\right)^{n},
\]
showing that $h_j=0$ if $j<0$, $h_0=1$, $h_1=e_1$, $h_2=e_1^2-e_2$,\ldots
Moreover the equality $H_r(t)E_r(t)=1$ implies the relation
\be
h_{n}+\sum_{1\leq i\leq r}(-1)^ie_ih_{n-i}=0,\qquad (n\in\ZZ),\label{eq:relhe}
\ee
which gives for   ($k\in\NN^*$) the well known relation  $h_k=\det(e_{j-i+1})_{1\leq i,j\leq k}$, 
under the convention $e_0=1$ and $e_j=0$ if  $j>r$,  \cite{MacDonald}.  The first few terms of $H_r$ in terms of the $x_i$s are (Cf.~\cite{KR} where the $h_i$s are called $S_i$)
\[
h_1=x_1, \quad h_2={x_1^2\over 2}+x_2,\quad h_3={x_1^3\over 3!}+x_1x_2+ x_3,\ldots
\]
 Let $\Ical_{H_r}$ be the ideal of $B_r$ generated by the relations~(\ref{eq:relhe}). Then
\[
B_r=\QQ[H_r]={\QQ[H_\infty]\over {\Ical}_{H_r}}\cong \QQ[h_1,\ldots,h_r].
\]
 In other words, if $r<\infty$,   each $h_{r+1+j}$ is an explicit polynomial expression in $(h_1,\ldots, h_r)$.
 
Similarly   the equality $E_r(t)^{-1}=\exp (X_r(t))$ implies 
\be
(r+1+j)x_{r+1+j}-(r+j)e_1x_{r+j}+\ldots+(-1)^rje_rx_{j+1}=0, \quad (j\geq 0)\label{eq;relationx}
\ee
 i.e. each $x_{r+1+j}$ can be expressed as a $B_r$-linear combination of $x_{j+1},\ldots, x_{r+j}$. By induction each $x_{r+1+j}$ is a polynomial in $x_1,\ldots, x_r$ with $\QQ$-coefficients. If $\Ical_{X_r}$ is the ideal of the relations~(\ref{eq;relationx}), then
\[
B_r=\QQ[X_r]={\QQ[X_\infty]\over {\Ical}_{X_r}}\cong \QQ[x_1,\ldots,x_r].
\]
Writing $h_n,x_n\in B_r$ means, respectively, the unique polynomial $h_n\,\,\mod \,\,\Ical_{H_r}\in \QQ[h_1,\ldots, h_r]$ and the unique polynomial
$x_n\,\,\mod\,\, \Ical_{X_r}\in \QQ[x_1,\ldots, x_r]$, where $x_n, h_n$ have been regarded as elements of $B_s$ with $s>n$.
For instance if $r=2$, then $h_3=h_1^3-2h_1h_2\in \QQ[h_1,h_2]$.
\claim{} 
If $\blamb=(\lambda_1,\ldots, \lambda_{r})\in \Pcal_r$  and  $H_r$ is as in~\ref{hr}, let
\[
\Delta_\blamb(H_r):=\det(h_{\lambda_{j}-j+i})_{1\leq i,j\leq r}\,\,\mod\,\, \Ical_{H_r}.
\]
%Then $(\Delta_{\blamb}(H_r))_{\blamb\in \Pcal_r}$ is a basis of the $\QQ$-vector space $\QQ[h_1,\ldots, h_r]\cong B_r=\QQ[H_\infty]/\Ical_{H_r}$.
\noindent
It is well known that \cite[p.~41]{MacDonald}:
\[
\QQ[h_1,\ldots, h_r]\cong B_r=\bigoplus_{\blamb\in \Pcal_r}\QQ\cdot \Delta_\blamb(H_r).
\]
%i.e. each polynomial in the $e_i$s can be written as a finite linear combination of $ \Delta_\blamb(H_r)$ with $\QQ$-coefficients. 

%\claim{\bf Remark.} If $r<\infty$,  let $\bfz_r:=(z_1,\ldots, z_r)$ be a set of indeterminates. Defines $(e_1(\bfz_r),\ldots, e_r(\bfz_r))$ via
%\[
%\sum_{i\geq 0}e_i(\bfz_r)t^i=\prod_{i=1}^r(1+z_it)
%\]
%then $e_j$ is the $j$th {\em elementary symmetric polynomial} in $z_1,\ldots, z_r$. In this case 
%\[
%h_j=\sum_{1\leq i_1\leq\ldots\leq i_r\leq r} z_{i_1}\ldots z_{i_r}
%\]  
%is the $j$th {\em complete symmetric polynomial} and 
%\[
%x_j/j=z_1^j+z_2^j+\ldots+z_r^j.
%\]
%the $j$th {\em Newton polynomial}, the sum of the $j$-th powers of $(z_1,\ldots, z_r)$.
\claim{} 
From now on let $A$ be any  $B_r$-algebra,  fixed once and for all. For $\phi=\sum_{n\geq 0} a_n t^n/n!$ define linear forms $U_k: A[[t]]\sra A$ by:
 \be
U_k(\phi)=a_{k} +\sum_{i\geq 1}(-1)^ie_ia_{k-i},\qquad (k\in\NN),\label{eq:U_k}
 \ee
 agreeing that $e_j=0$ if $j\geq r+1$ and $a_j=0$ if $j<0$.  
In particular $U_0(\phi)=a_0$, $U_1(\phi)=a_1-e_1a_0$, $U_2(\phi)=a_2-e_1a_1+e_2a_0$, \ldots
 Let
\be
u_j:=\sum_{n\geq 0}h_{n+j}{t^n\over n!}\in B_r[[t]], \quad (j\in\ZZ).\label{eq:uj}
\ee
One easily checks that $D^iu_j=u_{i+j}$ for $i\geq 0$ and $j\in\ZZ$.  By abuse of notation, we shall write $u_j\in A[[t]]$ instead of $u_j\otimes_{B_r}1_A\in A[[t]]$.
Notice that for $j\geq 0$
\be
u_{-j}={t^j\over j!}+h_1{t^{j+1}\over (j+1)!}+\ldots\label{eq:u-j}
\ee
and
\be
u_j=h_j+h_{j+1}t+h_{j+2}{t^2\over 2!}+\ldots
\ee
Obviously, all the $u_j$s are linearly independent over $\QQ$. Moreover:
\bclm{\bf Proposition.}\label{propdbas} {\em If $\phi\in A[[t]]$ then:
\be
\phi=\sum_{i\geq 0}U_i(\phi)u_{-i}.\label{eq:projphi}
\ee
}
\eclm
\proof Because of~(\ref{eq:u-j}), it is obvious that each $\phi\in A[[t]]$ can be written as an infinite linear combination $\phi=\sum_{i\geq 0} a_iu_{-i}$, for some $a_i\in A$. For each $j\geq 0$, the linearity of $U_j$ yields
$
U_j(\phi)=\sum_{i\geq 0} a_iU_j(u_{-i}).
$
By~(\ref{eq:u-j}) and Definition~(\ref{eq:U_k}),   $U_j(u_{-j})=1$.
If $j\neq -i$, instead
\[
U_j(u_{-i})=h_{n-i+j}+\sum_{1\leq k\leq r}(-1)^ie_kh_{n-i+j-k}=0,
\]
which is~(\ref{eq:relhe}) for $n-i+j$. Thus $U_j(\phi)=a_j$ and  each $\phi\in A[[t]]$ is the sum~(\ref{eq:projphi}) of its ``projections'' along  $u_{-i}$, for $i\in\NN$.\qed
%\claim{\bf Remark.} Relations~(\ref{eq:projphi}) can be  inverted. If $\phi$ is written as $\sum_{n\geq 0}a_nt^n/n!$, in fact,  it is easy to check that:
%\be
%a_k=U_k(\phi)+h_1U_{k-1}(\phi)+\ldots+ h_{k}U_0(\phi).\label{eq:invprojphi}
%\ee

The Proposition implies:
\bclm{\bf Corollary.} {\em For each $r\in \NN\cup\{\infty\}$ and each  $j\geq 0$:
\be
{t^j\over j!}=u_{-j}+\sum_{i\in [1,r]\cap \NN}(-1)^je_ju_{-j-i}.\label{eq:tjlcomb}
\ee
}
\eclm
\qed

\noindent
\begin{small}
For example, if $r=3$, one has $1=u_0-e_1u_{-1}+e_2u_{-2}+e_3u_{-3}$ and for $j\geq 1$
\[
\quad { t^{j}\over j!}=u_{-j}-e_1u_{-j-1}+e_2u_{-j-2}+e_3u_{-j-3}.
\]
\end{small}

\noindent
%\proof Just apply Proposition~\ref{propdbas}. 
Notice that if $r=\infty$ the r.h.s. of~(\ref{eq:tjlcomb})  is an infinite sum.

\claim{\bf Remark.} Relation~(\ref{eq:for1}) shows that each $\phi\in B_r[[t]]$ can be regarded as a  function of $(x_i)_{i\in [1,r]\cap \NN}$. 

%The proposition below will be useful later on.
%\bclm{\bf Lemma.}\label{lmma37} {\em If  $j\geq 0$ and $\phi\in B_r[[t]]$:
%\be
%{\partial U_j(\phi) \over \partial x_i}=U_j\left({\partial \phi\over \partial x_i}\right)+A_{j1}U_{j-1}(\phi)+\ldots+A_{j,j}U_{0}(\phi)\label{eq:dUjphi}
%\ee
%for some $A_{ji}\in B_r$ ($1\leq i\leq j$).
%}
%\eclm
%\proof 
%By definition of  $U_j$:
%\begin{eqnarray}
%{\partial U_j(\phi) \over \partial x_i}&=&{\partial\over \partial x_i}\left(a_j+\sum_{k\geq 1}(-1)^ke_ka_{j-k}\right)=\nonumber\\
%&=&{\partial a_j\over \partial x_i}+\sum_{k\geq 1}(-1)^ke_k{\partial a_{j-k} \over \partial x_i}+\sum_{k\geq 1}(-1)^k{\partial e_{j-k} \over \partial x_i}a_{j-k}=\nonumber \\
%&=&U_j\left({\partial \phi\over \partial x_i}\right)+\sum_{k\geq 1}(-1)^k{\partial e_k\over \partial x_i}a_{j-k}.\label{eq:dUjphip}
%\end{eqnarray}
%By writing each $a_{j-k}$  occurring in~(\ref{eq:dUjphip}) according to~(\ref{eq:invprojphi}) one obtains expression~(\ref{eq:dUjphi}).
%\qed

\section{Generic linear ODEs} \label{GODE}

\claim{} Let $t,z$ be  indeterminates over $B_r$ and  $e^{Dz}:=\exp(Dz)\in B_r[D][[z]]$. Clearly $e^{Dz}:B_r[[t]]\sra B_r[[t,z]]$ and
\be
U_r(e^{Dz}) =D^r-e_1D^{r-1}+\ldots+(-1)^re_r\in B_r[D]
\ee
is the {\em generic linear ordinary differential operator} of order $r$. Denote  $K_{r}:=\ker U_r(e^{Dz})$,   the submodule of $B_r[[t]]$ whose elements are  solutions of  $U_r(e^{Dz})y=0$,  the {\em generic linear homogeneous ODE}  of order $r$.

\bclm{\bf Proposition.}\label{propsol} {\em We have  $\phi\in K_r\otimes_{B_r}A$ $\iff$  $U_{n+r}(\phi)=0$ for all $n\geq 0$.
}
\eclm
\proof 
Write
\begin{eqnarray*}
U_{r}(e^{Dz})\phi&=&\hskip-2pt(D^r-e_1D^{r-1}+\ldots+(-1)^re_r)\sum_{n\geq 0}a_n{t^n\over n!}=\\
&=&\sum_{n\geq 0}\left(a_{n+r}-e_1a_{n+r-1}+\ldots+(-1)e_ra_{n}\right){t^n\over n!}.
\end{eqnarray*}
Thus $\phi\in K_r\otimes_{B_r}A$ if and only if $U_{n+r}(\phi)=0$ for all $n\geq 0$.\qed

\bclm{\bf Proposition.}\label{ubas} {\em If $r<\infty$ the $r$-tuple  $(u_0,u_{-1},\ldots , u_{-r+1})$ is a $B_r$-basis of $K_{r}$.}
\eclm
\proof
Clearly $(u_0,u_1,\ldots, u_{-r+1})$ are $B_r$-linear independent by~(\ref{eq:u-j}). They belong to $K_{r}$ because
$
U_{r+n}(u_{-j})=0
$
(Proposition~\ref{propdbas}),
 and are solutions to the generic linear ODE by Proposition~\ref{propsol}.
Then each $\phi\in K_r\subseteq A[[t]]$  can be written as 
\be
\phi=U_0(\phi)u_0+U_1(\phi)u_{-1}+\ldots+U_{r-1}(\phi)u_{-r+1},\label{eq:univcauchy}
\ee
applying once again Proposition~\ref{propdbas} and  Proposition~\ref{propsol}.
\qed

\noindent
Formula~(\ref{eq:univcauchy}) is called in~\cite{GS1} {\em Universal Cauchy formula}.

\claim{\bf Remark.} The solutions $(u_{-j})_{j\in [0,r]\cap\NN}$ of the generic linear ODE of order $r$ are {\em universal} in the following sense. For any associative commutative $\QQ$-algebra $A$ and each monic $P\in A[t]$ of degree $r$, then 
$
\ker P(D)=K_{r}\otimes_{B_r}A
$,
where $A$ is regarded  as a $B_r$-algebra via the unique morphism that maps $e_i\mapsto e_i(P)$. 

\claim{\bf Remark.} \label{lodeinfty} If $r=\infty$ we say that $(u_0,u_{-1},\ldots)$ is a fundamental system of solutions to the linear ODE of infinite order having $E_\infty:=(e_1,e_2,\ldots)$ as sequence of coefficients. In this sense, each element $\phi\in B_\infty[[t]]$ is a solution to the linear ODE of infinite order. The relation $\sum_{n\in\ZZ} h_nt^n=\exp(\sum x_it^i)$ easily implies that
\be
{\partial^i h_n\over \partial x^j_1}={\partial h_n\over \partial x_j}=h_{n-j}.\label{eq:derhni}
\ee
More generally:\[
{\partial u_j\over \partial x_i}=\sum_{n\geq 0}{\partial h_{n+j}\over \partial x_i}{t^n\over n!}=\sum_{n\geq 0}h_{n+j-i}{t^n\over n!}=u_{j-i}, \quad (j\in\ZZ).
\]
Thus for  $n\geq 1$,  $u_j$ is solution of the PDE
\[
\phi{\partial^4 \phi\over\partial x_n^4}-4{\partial \phi\over\partial x_n}{\partial \phi\over\partial x_n^3}+3\phi\left({\partial^2 \phi\over\partial x_n^2}\right)^2-3\phi \left({\partial^2 \phi\over\partial x_{2n}^2}\right)^2+3{\partial^2 \phi\over\partial x_{2n}^2}\phi+4{\partial \phi\over\partial x_{3n}}{\partial \phi\over\partial x_n}-4\phi{\partial^2 \phi\over\partial x_{3n}\partial x_n}=0,
\]
which is the bilinear form of the Kadomtsev-Petviashvilii (KP) equation
\be
{3\over 4}{\partial ^2 f\over \partial y^2}-{\partial \over \partial x}\left({\partial f\over \partial t}-{3\over 2}f{\partial f\over \partial x}-{1\over 4}{\partial^3 f\over \partial x^3}\right)=0,\label{eq:KP}
\ee
up to the substitution $x_n=x, x_{2n}=y, x_{3n}=t$ and putting $f(x,y,t)=2{\partial^2\over \partial x^2}(\log \phi)$. This is a nice way to phrase  the known fact that the complete symmetric polynomials in infinitely many indeterminates are solutions to the KP equation~(\ref{eq:KP}) (in the Hirota bilinear form-- see \cite[p.~75]{KR}). 
\section{Fermionic Fock spaces} \label{FFS}
\claim{} Given  $r\in \NN\cup\{\infty\}$, let $V_{r}:=\bigoplus_{j\in\ZZ}\QQ\cdot u_j$, with  $u_j\in B_{r}[[t]]$  as in~(\ref{eq:uj}). 
%For each $i\in\ZZ$ let 
%\[
%V_{r,i}:=\Span_\QQ[u_i, u_{i-1},u_{i-2},\ldots]
%\]
For each $i\in \ZZ$, let
\[
\Phi^r_{i}=u_{i}\w u_{i-1}\w\ldots\w u_{i-r+1}\w u_{i-r}\w u_{i-r-1}\w\ldots
\]
and for each $\blamb\in \Pcal_r$, let
\[
\Phi^r_{i,\blamb}:=u_{i+\lambda_{1}}\w u_{i-1+\lambda_{2}}\w\ldots\w u_{i-r+1+\lambda_{r}}\w \Phi^r_{-r}
\]
so that $\Phi^r_{i,0}=\Phi^r_{i}$. Let $F^r_{i}:=\bigoplus_{\blamb\in \Pcal_r}\QQ\cdot \Phi^r_{i,\blamb}$. Notice that $i\neq j$ implies that $F^r_{i}\cap F^r_{j}={\mathbf 0}$. The direct sum $\bigoplus_{i\in\ZZ}F^r_{i}$ is often denoted by  $\bw^{\infty/2}_\QQ V_r$ in the literature, and  is called {\em semi-infinite exterior power of $V_r$}. 
\bclm{\bf Remark.} If $r=\infty$, we write $\Phi_i$ instead of $\Phi^\infty_{ i}$ and $F_i$ instead of $F^\infty_{i}$. In particular 
 $F_{i}=\bigoplus_{\blamb\in \Pcal}\QQ\cdot \Phi_{i,\blamb}$, and a typical monomial of $F_{i}$ is
\[
u_{i+\lambda_{1}}\w\ldots\w u_{i-k+1+\lambda_{k}}\w \Phi_{i-k}
\]
where $\blamb=(\lambda_{1}\geq \lambda_{2}\geq\ldots\geq \lambda_{k})$ is a partition of length at most $k$, $k\geq 1$.
\eclm
It is a trivial remark that for  each $i$ and each $r$ there is a natural, although in principle not canonical, isomorphism $\sigma^r_i:F^r_{i}\sra B_r$ given by $\Phi^r_{i,\blamb}\mapsto \Delta_\blamb(H_r)$. Using that isomorphism each $F^r_{i}$ gets a structure of $B_r$-module of rank $1$ generated by $\Phi^r_{i}$. On the other hand $K_r=\ker U_r(e^{Dz})$ has a canonical structure of $B_r$-module of rank $r$ \cite[Theorem 2.1]{GS1}. It turns out that $\bw^rK_r$ itself is a free $B_r$-module of rank $1$ and thus $F^r_0:=\bw^rK_r\w \Phi^r_{-r}$ inherits a structure of rank $1$ free $B_r$-module. 
Our first remark is that such a structure coincides with that  induced by $\sigma^r_{0}$, via the Universal Cauchy Formula~(\ref{eq:univcauchy}).
\bclm{\bf Proposition (The Boson-Fermion Correspondence for $F^r_{0}$).} \label{pregmb}{\em The vector space $F^r_{0}$ has a canonical structure of free $B_{r}$-module of rank $1$ generated by $\Phi^r_{0}$.
}
\eclm
\proof
We shall construct a $B_r$-module morphism $F^r_{0}\sra B_r$ mapping $\Phi^r_{0}$ to $1$. We distinguish two cases. If $r<\infty$, one sees that 
\[
F^r_{0}=\bw^rK_{r}\w_\QQ\Phi^r_{-r}:=\Span_\QQ\{v_0\w v_1\w\ldots\w v_{r-1}\w \Phi^r_{-r}\,|\, v_i\in K_{r}\}.
\]
By the universal Cauchy formula~(\ref{eq:univcauchy}),  $v_{i}=\sum_{0\leq j\leq r-1}U_j(v_{i})u_{-j}$. For each typical monomial of $F^r_{0}$ one has
\[
 (v_0\w v_{1}\w \ldots\w v_{r-1}\w u_{-r}\w u_{-r-1}\w\ldots)\otimes_\QQ 1_{B_r}=
\]
\[
=\sum_{0\leq i_0\leq r-1}U_{i_0}(v_0)u_{-i_0}\w \sum_{0\leq i_{1}\leq r-1}U_{i_{1}}(v_1)u_{-i_1}\w\ldots \w\sum_{0\leq i_{r-1}\leq r-1}U_{i_{r-1}}(v_{r-1})u_{-i_{r-1}}\w\Phi^r_{-r}=
\]
\be
=\left| \matrix{U_{0}(v_0)&\ldots&U_0(v_{r-1})\cr
\vdots&\ddots&\vdots\cr
U_{r-1}(v_0)&\ldots&U_{r-1}(v_{r-1})
}
\right| u_0\w u_{-1}\w\ldots\w u_{-r+1}\w \Phi^r_{-r}=\det((U_{i}(v_{j})_{0\leq i,j\leq r-1})\Phi^r_{0},\label{eq:expbfc}
\ee
and $\det((U_{i}(v_{j})_{0\leq i,j\leq r-1})\in B_r$.  
Define $F^r_{0}\sra B_r$ by 
\[
v_0\w v_{1}\w \ldots\w v_{r-1}\w u_{-r}\w u_{-r-1}\w\ldots\mapsto {(v_0\w v_{1}\w \ldots\w v_{r-1}\w\Phi^r_{-r})\otimes_\QQ{1_{B_r}}\over \Phi^r_{0}}:=\det((U_{i}(v_{j})_{0\leq i,j\leq r-1})
\]
It is clearly an isomorphism. In fact  $\det((U_{i}(v_{j})_{0\leq i,j\leq r-1})=0$ implies $v_0,\ldots,$ $v_{r-1}$ are linearly dependent over $B_r$ and then $v_0\w\ldots\w v_{r-1}\w \Phi^r_{-r}=0$. Moreover $\Phi^r_{0}/\Phi^r_{0}=1$.

\noindent
If $r=\infty$ and $v_0,\ldots, v_{k-1}\in B_\infty$, then for each $0\leq i\leq k-1$
 \be
 v_i=\sum_{j\geq 0} U_j(v_i)u_{-j},\label{eq:infvuj}
 \ee
 which is in general an infinite linear combination. However substituting the expression~(\ref{eq:infvuj}) of $v_i$ into
$
v_0\w v_{1}\w \ldots\w v_{k-1}\w u_{-k}\w u_{-k-1}\w\ldots
 $
  one sees that all the summands  involving $u_{-j}$, with $j\geq k$, vanish due to skew-symmetry of the wedge product. Arguing as 
  in the  case $r<\infty$, one easily obtains the desired formula:
 
\medskip
  \hskip30pt $
v_0\w v_{1}\w \ldots\w v_{k}\w u_{-k-1}\w u_{-k-2}\w\ldots=\det((U_{i}(v_{j})_{0\leq i,j\leq k})\Phi_0.
  $
  \qed

\bclm{\bf Corollary.} {\em One has:
\[
{\Phi^r_{0,\blamb}\over \Phi^r_{0}}=\sigma^r_{0}(\Phi^r_{0,\blamb})=\Delta_\blamb(H_r).
\]
}
\eclm
\proof  If $r<\infty$ and  $v_i=u_{-i+\lambda_{i+1}}$, $0\leq i\leq r-1$ one has
\[
{\Phi^r_{0,\blamb}\over \Phi^r_{0}}=\det((U_{i}(u_{-j+\lambda_j})_{0\leq i,j\leq r-1}).
\]
By writing explicitly the determinant and using the definition of the $u_j$s and of the linear maps $U_i$ one  obtains:
\begin{eqnarray}
&&\left| \matrix{U_{0}(u_{0+\lambda_1})&\ldots&U_0(u_{-r+1+\lambda_{r}})\cr
\vdots&\ddots&\vdots\cr
U_{r-1}(u_{0+\lambda_1})&\ldots&U_{r-1}(u_{-r+1+\lambda_{r}})
}
\right|=\nonumber
\\ \nonumber\\
&=&\left|\matrix{h_{\lambda_1}&h_{\lambda_2-1}&\ldots&h_{\lambda_{r}-r+1}\cr
h_{\lambda_1+1}-e_1h_{\lambda_1}&h_{\lambda_2}-e_1h_{\lambda_2-1}&\ldots&h_{\lambda_{r}-r+2}-e_1h_{\lambda_{r}-r+1}\cr
\vdots&\vdots&\ddots&\vdots\cr
\sum_{i=0}^{r-1}(-1)^ie_ih_{\lambda_1+r-i}&\sum_{i=0}^{r-1}(-1)^ie_ih_{\lambda_2+r-1-i}&\ldots&\sum_{i=0}^{r-1}(-1)^ie_ih_{\lambda_{r}-i}
}\right|.\label{eq:dethuge}
\end{eqnarray}
Using the multilinearity and skew symmetry,  the determinant~(\ref{eq:dethuge}) simplifies into:
\[
\left| \matrix{h_{\lambda_1}&h_{\lambda_2-1}&\ldots&h_{\lambda_{r}-r+1}\cr
h_{\lambda_1+1}&h_{\lambda_2}&\ldots&h_{\lambda_{r}-r+2}\cr
\vdots&\vdots&\ddots&\vdots\cr
h_{\lambda_1+r-1}&h_{\lambda_2+r-2}&\ldots&h_{\lambda_{r}}\cr
}
\right|=\det(h_{\lambda_{j}-j+i})_{0\leq i,j\leq r-1}=\Delta_\blamb(H_{r}),
\]
as desired.  \qed

\claim{} \label{notationpieri} The module structure $B_r\otimes F^r_0\sra F^r_0$ induced by the Boson-Fermion correspondence is then defined by imposing the equalities
\[
\left\{\matrix{P\Phi^r_{0,\blamb}(H_r)&=&(P\cdot \Delta_\blamb(H_r))\Phi^r_0,\cr\cr \Delta_\blamb(H_r))\Phi^r_0&=&\Phi^r_{0,\blamb}.}\right.
\]
Since each $P\in B_r$ is a polynomial in $(e_i)_{i\in [1,r]\cap\NN}$ it suffices to know how to expand the product $e_i\Delta_\blamb(H_r)$ as a $\QQ$-linear combinations of Schur polynomials $\Delta_\blamb(H_r)$. This is prescribed by the following version of Pieri's formula:
\be
e_i\Delta_\blamb(H_r)=\Delta_{\blamb+i}(H_r)\label{eq:Pieri}
\ee
where for each integer $i\in [0,r]\cap \NN$ we set
\[
\Delta_{\blamb\pm i}(H_r)=\sum \Delta_{(\lambda_1\pm i_1,\ldots,\lambda_{r}\pm i_r)}(H_r)
\]
 the sum being  over  all $r$-tuples $(i_1,i_2,\ldots, i_r)$ such that $0\leq i_j\leq 1$, $\sum i_j=i$ and
\[
\lambda_1\pm i_1\geq \ldots\geq \lambda_{r}\pm i_r.
\]
(i.e. $\blamb\pm i$ is a partition).  
For example:
\[
\Delta_{(32)+2}(H_3)=\Delta_{(43)}(H_3)+\Delta_{(331)}(H_3)+\Delta_{(421)}(H_3)
\]
\[
\Delta_{(32)-1}(H_3)=\Delta_{(22)}(H_3)+\Delta_{(31)}(H_3)
\]
A minute of reflection shows that the action of $e_i$ can be described directly on $F^r_0$. Let 
\[
u_{i_0}\w u_{i_{-1}}\w \ldots \w u_{i_{-r+1}}\w \Phi^r_{-r}\in F^r_0,
\]
where
the indices $i_0,i_{-1},\ldots, i_{-r}$ are not necessarily in decreasing  order. Then:
\be
e_ju_{i_0}\w u_{i_{-1}}\w \ldots \w u_{i_{-r+1}}\w \Phi^r_{-r}=\sum_{(j_0,\ldots,j_{r-1})\in m(j)}u_{i_0+j_0}\w u_{i_{-1}+j_1}\w \ldots \w u_{i_{-r+1}+j_{r-1}}\w \Phi^r_{-r}\label{eq:periFro}
\ee
where $m(j):=\{(j_0,\ldots,j_{r-1})\in \NN^r\,|\, 0\leq j_i\leq 1, \sum j_i=j\}$. Clearly on the right hand side of~(\ref{eq:periFro}) some summands can vanish, whence the surviving partitions in formula~(\ref{eq:Pieri}).

 \claim{} Let $V_r^\vee=\hom_\QQ(V_r,\QQ)$. If $\alpha\in V^\vee_{r}$,   the  contraction map
$
 \alpha\in \hom_\QQ(F^r_{i+1},  F^r_{i})
 $ 
 is defined by:
 \[
 \alpha\lrcorner (u_{i_k}\w u_{i_{k-1}}\w\ldots)=\alpha(u_{i_k}) u_{i_{k-1}}\w u_{i_{k-2}}\w\ldots-\alpha(u_{i_{k-1}}) u_{i_{k}}\w u_{i_{k-2}}\w\ldots+\alpha(u_{i_{k-2}}) u_{i_{k}}\w u_{i_{k-1}}\ldots
 \]
 For each $i\in\ZZ$, let $u_i^\vee\in V^\vee_r$ defined by  $u_i^\vee(u_j)=\delta_{ij}$.
  Following \cite{KR}, define the formal Laurent series
  \[
  X(z)=\sum_{i\in\ZZ}u_iz^i\in V_{r}[[z,z^{-1}]]\qquad \mathrm{and} \qquad X^\vee(z)=\sum_{i\in\ZZ}u_i^\vee z^{-i}\in V_{r}^\vee[[z,z^{-1}]].
  \]
Consider
\[
\left\{\matrix{X(z)\w&:& F^r_{-1}&\lra& F^r_{0}[[z,z^{-1}]]&&\cr\cr
{}&&\Phi^r_{-1,\blamb}&\longmapsto&X(z)\w\Phi^r_{-1,\blamb}:=&&\hskip-14pt\sum_{i\in \ZZ} u_i\w\Phi^r_{-1,\blamb}\cdot z^i }\right.
\]
and 
\[
\left\{\matrix{X^\vee(z)\lrcorner &:& F^r_{1}&\lra& F^r_{0}[[z,z^{-1}]]\cr\cr
{}&&\Phi^r_{1,\blamb}&\longmapsto&X^\vee(z)\lrcorner \Phi^r_{1,\blamb}:=&&\hskip-14pt\sum_{i\in \ZZ} u_i\lrcorner\Phi^r_{1,\blamb}\cdot z^{-i}}\right.
\]
\claim{\bf Remark.} In the sequel we shall disregard the fermionic Fock spaces $F^r_i$ for $i\neq 0, 1,-1$, being irrelevant for the purposes of the present exposition. Most of the formulas deduced in the sequel, obtained for  $F^r_1, F^r_0$ and $F^r_{-1}$ only, can be easily generalized  for any triple $F^r_{i+1}, F^r_i$ and $F^r_{i-1}$ with no substantial change. This will be discussed in a forthcoming paper.
\claim{} The boson counterparts of the operatos $X(z)$ and $X^\vee(z)$ are the {\em vertex operators}
$
\Gamma_r(z): B_{r}\sra B_{r}[[z, z^{-1}]]$ and $\Gamma_r^\vee(z): B_{r}\sra B_{r}[[z^{-1},z]]
$
 defined as:
 \begin{eqnarray}
\Gamma_r(z) \Delta_\blamb(H_r) &=&{(X(z)\w \Phi^r_{-1,\blamb})\otimes_{\QQ}1_{B_r}\over \Phi^r_{0}},\label{eq:vvertex}\\ \nonumber \\
 \Gamma_r^\vee(z) \Delta_\blamb(H_r)&=&{(X^\vee(z)\lrcorner \Phi^r_{1,\blamb})\otimes_\QQ 1_{B_r}\over \Phi^r_{0}}.\label{eq:vvertexv}
\end{eqnarray}
 Notice that if one considers on $F^r_{1}$  the same $B_r$-module structure of $F^r_{0}$, it vanishes. This is why in~(\ref{eq:vvertexv}) the wedge is considered with respect to the $\QQ$-vector space structure.
The expression of $\Gamma_r(z)$ and $\Gamma_r^\vee(z)$ are very well known in the case when $r=\infty$: see~\cite{KR} and  Section~\ref{infty} where they are deduced in an alternative way. Next two sections are devoted to determine the shape and the properties of $\Gamma_r(z)$ and $\Gamma_r^\vee(z)$.

  \section{The Vertex Operator $\Gamma_r(z)$}\label{vertex}

\claim{} To describe the vertex operator $\Gamma_r(z):B_r\sra B_r[[z^{-1},z]]$
 we determine the action on each element of the distinguished basis $(\Delta_\blamb(H_r)\,|\, \blamb\in\Pcal_r)$ of $B_r$.  \bclm{\bf Lemma.} \label{fundlem}{\em For each $r\in \NN^*\cup\{\infty\}$:
\be
(X_r(z)\w\Phi^r_{-1,\blamb})\otimes_\QQ {1_{B_r}}:={1\over E_{r}(z)}\sum_{ i\in [0,r]\cap\NN}{E_i(z)\over z^i}u_{-i}\w \Phi^r_{-1,\blamb}.\label{eq:lemmafund}
\ee
}
\eclm
\proof
The Lemma will be proven  for $r<\infty$ in  a way that obviously extends to the case $r=\infty$.
Let $\bfu^r_{-1,\blamb}=u_{-1+\lambda_1}\w\ldots \w u_{-r+\lambda_r}$. Then
\be
X(z)\w \Phi^r_{-1,\blamb}=X(z)\w \bfu^r_{-1,\blamb}\w \Phi^r_{-r-1}=\sum_{i\geq -r}z^iu_i\w\bfu^r_{-1,\blamb}\w \Phi^r_{-r-1}\label{eq:mkssense}
\ee
If $r<\infty$, the summation index of~(\ref{eq:mkssense}) runs over all integers $\geq -r$ because $u_i\w  \Phi^r_{-r-1}=0$ for all $i\leq -r-1$.
Now:
\[
\sum_{i\geq -r}u_iz^i={u_{-r}\over z^r}+\sum_{0\leq i\leq r-1} {u_{-i}\over z^i}+\sum_{j\geq 1}u_jz^j={u_{-r}\over z^r}+\sum_{0\leq i\leq r-1} {u_{-i}\over z^i}+\sum_{j\geq 1}\sum_{0\leq i\leq r-1} U_i(u_j)u_{-i}z^j,
\]
where we wrote $u_j=\sum_{0\leq i\leq r-1}U_i(u_j)u_{-i}$ using Proposition~\ref{propdbas} and the fact that $U_i(u_j)=0$ if $i\geq r$, because $u_i$ is solution to the generic linear ODE of order $r$.
By suitably grouping the summands one obtains:
\be
\sum_{i\geq -r}u_iz^i={u_{-r}\over z^r}+\sum_{0\leq i\leq r-1}\left({1\over z^i}+\sum_{j\geq 1}U_i(u_j)z^j\right) u_{-i}.\label{eq:vttre}
\ee
Observes now that, since $U_{i+j}(u_0)=0$,
\[
U_i(u_j)=h_{i+j}-e_1h_{i+j-1}+\ldots+e_ih_j=-\sum_{i+1\leq k\leq r-1}(-1)^ke_kh_{j+i-k}
\]
and then~(\ref{eq:vttre}) can be written as 
\begin{eqnarray*}
\sum_{i\geq -r}u_iz^i&=&{u_{-r}\over z^r}+\sum_{0\leq i\leq r-1}\left({1\over z^i}-\sum_{j\geq 1}\sum_{i+1\leq k\leq r-1}(-1)^ke_kh_{j+i-k}z^j\right) u_{-i}=\\
&=&{u_{-r}\over z^r}+\sum_{0\leq i\leq r-1}\left({1\over z^i}-\sum_{i+1\leq k\leq r-1}(-1)^ke_k\sum_{j\geq 1}h_{j+i-k}z^j\right) u_{-i}=\\
&=&{u_{-r}\over z^r}+\sum_{0\leq i\leq r-1}\left({1\over z^i}-\sum_{i+1\leq k\leq r-1}(-1)^ke_k\sum_{j\geq 0}h_{j+1+i-k}z^{j+1}\right) u_{-i}=\\
&=&{u_{-r}\over z^r}+\sum_{0\leq i\leq r-1}\left({1\over z^i}-\sum_{i+1\leq k\leq r-1}(-1)^ke_k\sum_{j\geq 0}h_{j}z^{j+ k-i}\right) u_{-i}=\\
&=&{u_{-r}\over z^r}+\sum_{0\leq i\leq r-1}{1\over z^i}\left(1-\sum_{i+1\leq k\leq r-1}(-1)^ke_kz^k\sum_{j\geq 0}h_{j}z^{j}\right) u_{-i}=\\
&=&{u_{-r}\over z^r}+\sum_{0\leq i\leq r-1}{1\over z^i}\left({E_r(z)-\sum_{i+1\leq k\leq r-1}(-1)^ke_kz^k\over E_r(z)}\right) u_{-i}=
\end{eqnarray*}
\begin{eqnarray*}
&=&{u_{-r}\over z^r}+\sum_{0\leq i\leq r-1}{E_i(z)\over E_r(z)}{u_{-i}\over z^i}={1\over E_r(z)}\sum_{0\leq i\leq r} E_i(z){u_{-i}\over z^i}
\end{eqnarray*}
which proves~(\ref{eq:lemmafund}). If $r=\infty$ one has
\[
\sum_{i\geq -r}u_iz^i=\sum_{i\in\ZZ}u_iz^i=\sum_{i\geq 0} u_{-i}z^{-i}+\sum_{j\geq 0} u_{j}z^{j}
\]
and the same proof works in this case as well, up to expressing each $u_j$ as a (infinite) linear combination of the $u_{-i}$'s.\qed

\bclm{\bf Proposition.} \label{cortj} {\em Let $\blamb\in \Pcal_k$, $k\in [0,r]\cap \NN$. Then
\be
(X(z)\w \Phi^r_{-1,\blamb})\otimes_\QQ{1_{B_r}}={1\over E_r(z)}\cdot \exp\left({t\over z}\right)\w \Phi^r_{-1,\blamb}.\label{eq:expvert}
\ee

 }
 \eclm
 \proof
 By Lemma~\ref{fundlem}:
 \[
 (X(z)\w\Phi^r_{-1,\blamb})\otimes_\QQ{1_{B_r}}={1\over E_r(z)}\left( u_0+{u_{-1}\over z}E_1(z)+\ldots+{u_{-k}\over z^k}E_k(z)\right)\w \Phi^r_{-1,k}
 \]
  because $u_{-k-j}\w \Phi^r_{-1,\blamb}=0$ for all $j\geq 1$. But:
 \begin{eqnarray*}
  u_0+{u_{-1}\over z}E_1(z)+\ldots+{u_{-k}\over z^k}E_k(z)&=&\sum_{j=0}^k{1\over z_j}\left(u_{-j}+\sum_{i=1}^{k-j}(-1)^ie_iu_{-j-i}\right)=\\
  &=&u_0-e_1u_{-1}+\ldots+(-1)^ke_ku_{-k}+\\
  &+&{1\over z}\left(u_{-1}-e_1u_{-2}+\ldots+(-1)^{k-1}e_{k-1}u_{-k}\right)+\\
  &+&\ldots +\\
 &+&{1\over z^k}\left(u_{-k}-e_1u_{-k-1}+\ldots+(-1)^{k-1}e_{k-1}u_{-k}\right).
\end{eqnarray*}
 Now, for each $0\leq j\leq k$, because of~(\ref{eq:tjlcomb}):
 \[
u_{-j}+\sum_{i=1}^{k-j}(-1)^ie_iu_{-j-i}={t^j\over j!}+\sum_{p\geq 1}(-)^pe_{k-j+p}u_{-k-p}
 \]
 and so
 \[
 u_{-j}+\sum_{i=1}^{k-j}(-1)^ie_iu_{-j-i}\w \Phi^r_{-1,\blamb}=\left({t^j\over j!}+\sum_{p\geq 1}(-)^pe_{k-j+p}u_{-k-p}\right)\w\Phi^r_{-1,\blamb}={t^j\over j!}\w \Phi^r_{-1,\blamb}.
 \]
 Notice now that
 $
\displaystyle{ t^{k+1}\over (k+1)!}\w\Phi_{-1,\blamb}^r=0
 $
 because $t^{k+1}/(k+1)!$ is a linear combination of $u_{-k-j}$, $j\geq 1$. This prove that formula~(\ref{eq:lemmafund}) can be put in the form~(\ref{eq:expvert}), and the claim is proven. \qed

\claim{} For convenience we define
 $
 G_r(z): B_r\sra B_r[z^{-1}, z]]
 $
 through the equality 
 \[
 G_r(z)\Delta_\blamb(H_r)=E_r(z)(\Gamma_r(z)\Delta_\blamb(H_r)).
 \]
Lemma~\ref{fundlem} says that indeed $G_r(z)$ takes values in the polynomial ring $B_r[z^{-1}]$.
 
 \bclm{\bf Theorem.}\label{thmGrz} {\em
Notation as in~\ref{notationpieri}.  For each partition $\blamb$ of length $k\in [0,r]\cap\NN$:
 \be
 G_r(z)\Delta_\blamb(H_r)=\Delta_{\blamb}(H_r)-{1\over z}\Delta_{\blamb-1}(H_r)+{1\over z^2}\Delta_{\blamb-2}(H_r)+\ldots+(-1)^k{1\over z^k}\Delta_{\blamb-k}(H_r).\label{eq:forvertexg}
 \ee
 }
 \eclm
 \proof
 \medskip
\noindent
If  $\blamb=(\lambda_1,\ldots, \lambda_k)$ is a partition of length $k\leq r$, then:
\be
G_r(z)\Delta_\blamb(H_r)\Phi^r_{0}=\exp\left({t\over z}\right)\w \Phi^r_{-1,\blamb}+ =1\w \Phi^r_{-1,\blamb}+ {t\over z}\w \Phi^r_{-1,\blamb}+\ldots+ {1\over z^k}{t^k\over k!} \w\Phi^r_{-1,\blamb}\label{eq:prvexpvrt}
\ee
Now, for each $0\leq j\leq k$:
\begin{eqnarray*}
{t^j\over j!}\w \Phi^r_{-1,\blamb}&=&\left(u_{-j}-e_1u_{-j-1}+\ldots+(-1)^{k-j}e_{k-j}u_{-k}\right)\w\Phi^r_{-1,\blamb}=\\ 
&=&\left(u_{-j}\w\Phi^r_{-1,\blamb}-e_1u_{-j-1}\w\Phi^r_{-1,\blamb}+ \ldots+(-1)^{k-j}e_{k-j}u_{-k}\w\Phi^r_{-1,\blamb}\right)=\\ \\
&=&u_{-j}\w\Phi^r_{-1,\blamb}-(u_{-j}\w\Phi^r_{-1,\blamb}+u_{-j-1}\w \Phi^r_{-1,\blamb+1})+\\
&+&(u_{-j-1}\w \Phi^r_{-1,\blamb+1}+u_{-j-2}\w \Phi^r_{-1,\blamb+2})+\\
&+&\ldots+\\
&+&(-1)^{k-j}(u_{-k+1}\w\Phi^r_{-1,\blamb+k-j-1}+u_{-k}\w\Phi^r_{-1,\blamb+k-j})=\\
\end{eqnarray*}
so that only the term
$
(-1)^{k-j}u_{-k}\w\Phi^r_{-1,\blamb+k-j}
$
survives to cancelation. This last term can be written as follows:
\[
(-1)^{k-j}u_{-k}\w\Phi^r_{-1,\blamb+k-j}=(-1)^{k-j}u_{-k}\w\sum_{(j_1,\ldots, j_k)\in m(j)} u_{\lambda_1+j_1}\w\ldots\w u_{-k+1\lambda_{k}+j_k}\w \Phi^r_{-k-1}=
\]
\[
=(-1)^{j}\sum_{(j_1,\ldots, j_k)\in m(j)} u_{\lambda_1+j_1}\w u_{-1+\lambda_2+j_2}\w\ldots\w u_{-k+1+\lambda_{k}+j_k}\w u_{-k}\w \Phi^r_{-k-1}=
\]
\[
=(-1)^j\sum_{(j_1,\ldots, j_k)\in m(j)}\Delta_{(\lambda_1+j_1,\ldots,\lambda_k+j_k)}(H_r)\Phi^r_{0,\blamb}=(-1)^j\Delta_{\blamb-i}(H_r)\Phi^r_{0}.
\]
Substitution into~(\ref{eq:prvexpvrt}) gives~(\ref{eq:forvertexg}).\qed

\noindent
Applying Theorem~\ref{thmGrz} to the case $k=1$ one obtains:

\bclm{\bf Corollary.}\label{corhnz} {\em For each $r\geq 1$
\be
\Gamma_r(z)h_n={1\over E_r(z)}\left(h_n-{h_{n-1}\over z}\right).\label{eq:hnhnmunoz}
\ee
i.e.
\[
G_r(z)h_n=h_n-{h_{n-1}\over z}
\]
\qed
}
\eclm
Let $G_r(z)H_r$ be the sequence $(1,G_r(z)h_1, G_r(z)h_2,\ldots)$. Using Corollary~\ref{corhnz} it is easily checked that
\[
\Delta_\blamb(G_r(z)H_r)=\Delta_\blamb(H_r)-{1\over z}\Delta_{\blamb -1}(H_r)+\ldots+(-1)^r{1\over z^r}\Delta_{\blamb -r}(H_r)
\]
which so proves the first of our main results.

\bclm{\bf Theorem.}\label{mnth1} {\em The operator $G_r(z)$ commutes with taking $\Delta_\blamb$:
\be
G_r(z)\Delta_\blamb(H_r)=\Delta_\blamb(G_r(z)H_r)\label{eq:Gamrtrmk}
\ee
and then:
 \be
\Gamma_r(z)\Delta_\blamb(H_r)={1\over E_r(z)}{\Delta_\blamb(G_r(z)H_r)}.\label{eq:Gamcompart}
 \ee
 \qed
}
\eclm

\claim{\bf Remark.}\label{rmkdetG} The vector space  $V_r=\Span_\QQ(u_i)_{i\in\ZZ}$ is naturally a  $B_r$-module via the multiplication
\[
Pu_j:=P\sum_{n\geq 0}h_{n+j}{t^n\over n!}=\sum_{n\geq 0} Ph_n{t^n\over n!},\qquad (P\in B_r).
\]
One may so define $\widetilde{G}_r(z):V_r\sra V_r[z^{-1}]$ as
\[
\widetilde{G}_r(z)u_j=\sum_{n\geq 0} G_r(z)h_n{t^n\over n!}=u_j-{1\over z}\,u_{j-1}
\]
\noindent
Then~(\ref{eq:Gamrtrmk}) says that
\begin{eqnarray}
(G_r(z)\Delta_\blamb(H_r))\Phi^r_0&=&(\det \widetilde{G}_r(z))\cdot \Phi^r_{0,\blamb}=\nonumber\\&=&\widetilde{G}_r(z)u_{\lambda_1}\w \widetilde{G}_r(z) u_{-1+\lambda_2}\w\ldots\w\widetilde{G}_r (z)u_{-r+1+\lambda_r}\w \Phi^{r}_{-r}.
\end{eqnarray}

\bclm{\bf Corollary.}\label{qucom} {\em Let $h_{i_1}\cdot \ldots\cdot h_{i_s}$ be an arbitrary product of terms of $H_r$, with $s\leq r$. Then
\[
G_r(z)(h_{i_1}\cdot\ldots\cdot h_{i_s})=G_r(z) h_{i_1}\cdot\ldots\cdot G_r(z)h_{i_s}
\] 
}
\eclm
\proof Since $s\leq r$, each monomial in the $h_j$ is a $\ZZ$-linear combination of Schur polynomials associated to partitions of length at most $r$.
Suppose
$
h_{i_1}\cdot\ldots\cdot h_{i_s}=\sum_{|\blamb|\leq s}a_\blamb\Delta_\blamb(H_r)
$. 
Then
\begin{eqnarray*}
\hskip105pt G_r(z)(h_{i_1}\cdot\ldots\cdot h_{i_s})&=&\sum_{|\blamb|\leq s}a_\blamb G_r(z)\Delta_\blamb(H_r)=\\
&=&\sum_{|\blamb|\leq s}a_\blamb\Delta_\blamb(G_r(z)H_r)=\\&=&G_r(z)h_{i_1}\cdot\ldots\cdot G_r(z)h_{i_r}.\hskip105pt \qed
\end{eqnarray*}

 \claim{\bf Example.} Notice  that if  $r<\infty$, $G_r(z)$ is not a ring homomorphism. Consider e.g.\linebreak $\Gamma_1(z) :B_1\sra B_1[[z, z^{-1}]$. Then
 \[
 \Gamma_1(z) h_2={X_1(z)\w u_1\w u_{-2}\w\ldots\over \Phi^1_{0}}={-u_{1}z\w u_{-1}\w u_{-2}\over\Phi^1_{0}}={h_1\over z}={1\over E_1(z)}\left (h_2-{h_1\over z}\right)
 \]
 Similarly
 \[
  \Gamma_1(z)h_1={X_1(z)\w u_0\w u_{-2}\w\ldots\over \Phi^1_{0}}={-(1/z)u_{0}\w {u_{-1}}\w u_{-2}\over\Phi^1_{0}}=-{1\over z}={1\over E_1(z)}\left (h_1-{1\over z}\right)
 \]
 Therefore
 \[
 G_1(z)h_2=h_2-{h_1\over z}\qquad \mathrm{and}\qquad G_1(z)h_1=h_1-{1\over z}
 \]
 However in $B_1$  the relation $h_2=h_1^2$ holds. Hence
 \[
 h_2-{h_1\over z}=G_1(z)(h_1^2)\neq G(z)h_1\cdot G(z)h_1=\left(h_1-{1\over z}\right)^2=h_2-{2h_1\over z}+{1\over z^2}.
 \]

\bclm{\bf Corollary.}\label{liminf}{ \em If $r=\infty$, then $G_\infty(z):B_\infty\sra B_\infty[z^{-1}]$ is a ring homomorphism.\qed
}
\eclm
\noindent

\section{The Vertex Operator $\Gamma_r^\vee(z)$}\label{secgammavee}

 Let $G_r^\vee(z):B_r\sra B_r[[z^{-1},z]]$ be the $\QQ$-homomorphism
\be
\Delta_\blamb(H_r)\mapsto G_r^\vee(z)\Delta_\blamb(H_r):={z(\Gamma_r^\vee(z)\Delta_\blamb(H_r))\over E_r(z)}=[{z\Gamma_r^\vee(z)\Delta_\blamb(H_r)}]\cdot \sum_{n\geq 0}h_nz^n,
\ee
where $\Gamma^\vee_r(z):B_r\sra B_r[[z^{-1},z]]$ is as in~(\ref{eq:vvertexv}).

\bclm{\bf Lemma.} \label{scndlm} { \em For each $n\geq 0$
\be
G^\vee_r(z) h_n=\sum_{i\geq 0}{h_{n-i}\over z^i}=h_n+{h_{n-1}\over z}+\ldots+{h_1\over z^{n-1}}+{1\over z^n}.\label{eq:grveezhn}
\ee
}
\eclm
\proof
By definition of $G_r^\vee(z)$ one has
\[
E_r(z)(G^\vee_r(z)h_n)\Phi^{r}_0=z(\Gamma_r(z)h_n)\Phi^r_0=zX^\vee_r(z) \lrcorner h_n\Phi^r_1.
\]
For $r<\infty$:
\begin{eqnarray*}
zX_r^\vee(z)\lrcorner (h_n u_{1}\w u_{0}\w\ldots\w u_{-r+1}\w \Phi^{r}_{-r})&=&zX^\vee_r(z)\lrcorner (u_{1+n}\w u_0\w u_{-1}\ldots\w u_{-r+1}\w \Phi^r_{-r})=\\
&=&z^{-n}\Phi_0^{r}-u_{1+n}\w u_{-1}\w u_{-2}\w u_{-3}\ldots+\\
&+& z\cdot u_{1+n}\w u_{0}\w u_{-2}\w\ldots+\ldots+\\
&+&\ldots+\\
&+&z^{r}u_{1+n}\w u_{0}\w u_{-1}\w\ldots\w u_{r-2}\w \Phi^r_{-r}=
\end{eqnarray*}
\[
=(z^{-n}-h_{n+1}+\Delta_{(n+1,1)}(H_r)z^2+\ldots+(-1)^r\Delta_{(n+1,1^{r-1})}(H_r)z^{r})\Phi^r_0
\]
There are no further terms, because $\Delta_\blamb(H_r)=0$ if $\ell(\blamb)>r$.
An easy check shows that $(-1)^j\Delta_{(n+1, 1^j)}(H_r)=U_j(u_{n+1})$, so that
\be
zX^\vee(z)\lrcorner h_n\Phi^r_{1}=(z^{-n}-U_0(u_{n+1})z-U_1(u_{n+1})z^2+\ldots -U_{r-1}(u_{n+1})z^{r})\Phi^r_0\label{eq:sumxcheck}
\ee
On the other hand
\begin{eqnarray*}
&{}&E_r(z)\left(h_n+{h_{n-1}\over z}+\ldots+{h_1\over z^{n-1}}+{1\over z^n}\right)=\\
&=&z^{-n}E_r(z)(1+h_1z+h_2z^2+\ldots+h_nz^n)=\\
&=&z^{-n}E_r(z)\left({1\over E_r(z)}-\sum_{p\geq n+1}h_pz^p\right)=\\
&=&z^{-n}E_r(z)\left({1\over E_r(z)}-z^{n+1}{U_0(u_{n+1})+U_1(u_{n+1})z+\ldots+U_{r-1}(u_{n+1})z^{r-1}\over E_r(z)}\right)=\\
&=&z^{-n}-U_0(u_{n+1})z- U_1(u_{n+1})z^2-\ldots -U_{r-1}(u_{n+1})z^{r}={zX^\vee(z)\lrcorner h_n\Phi^r_{1}\over \Phi^r_0}
\end{eqnarray*}
and this proves~(\ref{eq:grveezhn}). \qed

\noindent
More generally:
\bclm{\bf Lemma.}\label{72} {\em  Let $(\lambda_1,\ldots, \lambda_r)$ be any partition of length at most $r$. Then:
\be
z\Gamma_r^\vee(z)\Delta_{(\lambda_1,\ldots,\lambda_r)}(H_r)=\left|\matrix{z^{-\lambda_1}&z^{1-\lambda_2}&\ldots&z^{r-1-\lambda_r}\cr
h_{\lambda_1+1}&h_{\lambda_2}&\ldots&h_{\lambda_r+r-2}\cr
\vdots&\vdots&\ddots&\vdots\cr
h_{\lambda_1+r-1}&h_{\lambda_2+r-2}&\ldots&h_{\lambda_r}}
\right|+\Delta_{(\lambda_1+1,\ldots, \lambda_r+1)}(H_r).\label{eq:detgavee}
\ee

%\begin{eqnarray*}
%\Gamma_r^\vee(z)\Delta_{(\lambda_1,\ldots,\lambda_r)}(H_r)&=&z^{-\lambda_1}\Delta_{(\lambda_2,\ldots,\lambda_r)}(H_r)-z^{1-\lambda_2}\Delta_{(\lambda_1+1,\lambda_3,\ldots,\lambda_r)}(H_r)+\\
%&+&\ldots+\\
%&+&(-1)^{r-1}z^{r-1+\lambda_{r}}\Delta_{\lambda_1,\ldots,\lambda_{r-1}}(H_r)-(-1)^{r-1}z^r\Delta_{(\lambda_1-1.\ldots,\lambda_r-1)}(H_r)=\\
%&=&\sum_{j=1}^r(-1)^{j-1}z^{j-1+\lambda_j}\Delta_{(\lambda_1+1,\ldots,\lambda_{j-1}+1,\lambda_{j+1},\ldots, \lambda_r)}(H_r)+(-1)^{r}z^r(-1)^{r-1}z^r\Delta_{(\lambda_1-1.\ldots,\lambda_r-1)}(H_r)
%\end{eqnarray*}
%
where as usual one sets $h_j=0$ if $j<0$.} 
\eclm

\proof
The proof of the equality is straightforward, as it  merely consists in expanding  the definition~(\ref{eq:vvertexv}) of $\Gamma^\vee_r(z)$.

\begin{eqnarray*}
zX_r^\vee(z)\lrcorner \Delta_\blamb(H_r)\Phi^r_{1,0}&=&zX^\vee(z)\lrcorner (u_{1+\lambda_1}\w u_{\lambda_2}\w\ldots \w u_{-r+1+\lambda_r}\w u_{-r}\w \Phi^r_{1,-r-1})=\\
&=&z^{-\lambda_1}\Delta_{(\lambda_2,\ldots,\lambda_r)}(H_r)-z^{1-\lambda_2}\Delta_{(\lambda_1+1,\lambda_3,\ldots,\lambda_r)}(H_r)+\\
&+&\ldots+\\
&+&(-1)^{r-1}z^{r-1+\lambda_{r}}\Delta_{\lambda_1,\ldots,\lambda_{r-1}}(H_r)-(-1)^{r-1}z^r\Delta_{(\lambda_1-1.\ldots,\lambda_r-1)}(H_r)=
\end{eqnarray*}
\be
=\sum_{j=1}^r(-1)^{j-1}z^{j-1+\lambda_j}\Delta_{(\lambda_1+1,\ldots,\lambda_{j-1}+1,\lambda_{j+1},\ldots, \lambda_r)}(H_r)+(-1)^{r}z^r(-1)^{r-1}z^r\Delta_{(\lambda_1-1.\ldots,\lambda_r-1)}(H_r),\label{eq:prvgvee}
\ee
and the first summand of~(\ref{eq:prvgvee}) is precisely the determinant occurring in~(\ref{eq:detgavee}).\qed

\bclm{\bf Theorem.} \label{mnthm2} {\em The operator $G_r^\vee(z)$ commutes with taking Schur determinats, i.e.:
\[
G^\vee_r(z)\cdot \Delta_\blamb(H_r)=\Delta_\blamb(G^\vee_r(z)H_r).
\]
Therefore
\[
\Gamma_r^\vee(z)\Delta_\blamb(H_r)={1\over z}E_r(z)\cdot \Delta_\blamb(G_r(z)H_r).
\]
}
\eclm

\proof
Again by definition of $G_r^\vee(z)$ one has
\[
G^\vee_r(z)\Delta_\blamb(H_r)={z\Gamma_r(z)\Delta_\blamb(z)\over E_r(z)}={z\Gamma_r(z)\Delta_\blamb(H_r)}\sum_{n\geq 0}h_nz^n.
\]
Using Proposition~\ref{72}:
\begin{eqnarray}
\hskip-13pt G_r^\vee(z)\Delta_\blamb(H_r)\hskip-7pt&=&\hskip-7pt z\Gamma_r^\vee(z)\Delta_{(\lambda_1,\ldots,\lambda_r)}(H_r)=\nonumber\\
 &=&\hskip-12pt\left(\left|\matrix{z^{-\lambda_1}&\hskip-10pt z^{1-\lambda_2}&\hskip-10pt \ldots&\hskip-10pt z^{r-1-\lambda_r}\cr
h_{\lambda_1+1}&\hskip-10pt h_{\lambda_2}&\hskip-10pt \ldots&\hskip-10pt h_{\lambda_r+r-2}\cr
\vdots&\vdots&\ddots&\vdots\cr
h_{\lambda_1+r-1}&\hskip-10pt h_{\lambda_2+r-2}&\hskip-10pt\ldots&\hskip-10pt h_{\lambda_r}}
\right|\hskip-3pt+\hskip-3pt(\hskip-1pt-1)^r\hskip-1pt\Delta_{(\lambda_1+1,\ldots, \lambda_r+1)}(H_r)z^r\hskip-3pt \right)\hskip-3pt \cdot\hskip-3pt\sum_{n\geq 0}h_nz^n.\label{eq:sgnr}
\end{eqnarray}
The  key computational remark is:  the coefficient of $z^n$, $n\in\ZZ$, in expression~(\ref{eq:sgnr}) is
\[
H(\blamb+ n)+(-1)^rh_{n-r}\Delta_{(\lambda_1+1,\ldots, \lambda_r+1)}(H_r)
\]
where we set, for sake of brevity
\be
H(\blamb+n):=\left|\matrix{h_{\lambda_i+1-i+n}\cr h_{\lambda_i+2-i}\cr\vdots\cr h_{\lambda_i+r-i}}
\right|_{1\leq i\leq r}=\left|\matrix{h_{\lambda_1+n}&h_{\lambda_2-1+n}&\ldots&h_{\lambda_r-r+1+n}\cr
h_{\lambda_1+1}&h_{\lambda_2}&\ldots&h_{\lambda_r+r-2}\cr
\vdots&\vdots&\ddots&\vdots\cr
h_{\lambda_1+r-1}&h_{\lambda_2+r-2}&\ldots&h_{\lambda_r}}
\right|\label{eq:skbrev}
\ee
and  $h_j=0$ if $j<0$. We claim that for  $n>0$
\[
H(\blamb+n)+(-1)^rh_{n-r}\Delta_{(\lambda_1+1,\ldots, \lambda_r+1)}(H_r)=0.
\]
  In fact if $1\leq n\leq r-1$, one has  $h_{n-r}\Delta_{(\lambda_1+1,\ldots, \lambda_r+1)}=0$, since $h_{n-r}=0$, while the first row of $H(\blamb+n)$ is equal to its $(n+1)$-th row  and so vanishes by skew-symmetry. For $n=r$ one has
 \[
 H(\blamb+r)+ (-1)^r\Delta_{(\lambda_1+1,\ldots, \lambda_r+1)}(H_r)=0
 \]
 as an immediate check shows (substitute $n=r$ in~(\ref{eq:skbrev})). For  $1\leq n-r\leq r-1$ one has:
\[
h_{\lambda_i+1-i+n}= \sum_{j=1}^{n-r}(-1)^{j-1}e_jh_{\lambda_i-i+1-j+n}+(-1)^{n-r+1} e_{n-r+1}h_{\lambda_i-i+r}+\ldots+ (-1)^re_rh_{\lambda_i+1-i+n-r}
\]
i.e.
\[
H(\blamb+n)=\left|\matrix{ \sum_{j=1}^{n-r}(-1)^{j-1}e_jh_{\lambda_i-i+1-j+n}+(-1)^{n-r+1} e_{n-r+1}h_{\lambda_i-i+r}+\ldots+ (-1)^re_rh_{\lambda_i+1-i+n-r}\cr  h_{\lambda_i+2-i}\cr\vdots\cr h_{\lambda_i+r-i}} \right|=
\]
\be
=\left|\matrix{ \sum_{j=1}^{n-r}(-1)^{j-1}e_jh_{\lambda_i-i+1-j+n}\cr  h_{\lambda_i+2-i}\cr\vdots\cr h_{\lambda_i+r-i}} \right|+ \left|\matrix{ (-1)^{n-r+1} e_{n-r+1}h_{\lambda_i-i+r}+\ldots+ (-1)^re_rh_{\lambda_i+1-i+n-r}\cr  h_{\lambda_i+2-i}\cr\vdots\cr h_{\lambda_i+r-i}} \right|\label{eq:scndsum}
\ee
The second summand in~(\ref{eq:scndsum}) vanishes because linearity and skewsymmetry of the determinant. Thus:
\[
H(\blamb+n)=e_1H(\blamb-1+n)+\ldots-(-1)^{n-r}e_{n-r}H(\blamb-r+n).
\]
In particular
\[
H(\blamb+r+1)+(-1)^rh_1\Delta_{(\lambda_1+1,\ldots.\lambda_r+1)}(H_r)=e_1H(\blamb+r)+(-1)^rh_1\Delta_{(\lambda_1+1,\ldots.\lambda_r+1)}(H_r)=0.
\]
By induction, for all $1\leq n-r\leq r-1$:
\[
H(\blamb+n)+(-1)^rh_{n-r}\Delta_{(\lambda_1+1,\ldots.\lambda_r+1)}(H_r)=
\]
\[
=\sum_{j=1}^{n-r}(-1)^{j-1}e_jH(\blamb+n-j)-\sum_{j=1}^{n-r}(-1)^{j-1}e_jh_{n-r-j}\Delta_{(\lambda_1+1,\ldots.\lambda_r+1)}(H_r)=
\]
\[
=\sum_{j=1}^{n-r}(-1)e_j\left(H(\blamb+n-j)-h_{n-r-j}\Delta_{(\lambda_1+1,\ldots.\lambda_r+1)}(H_r)\right)=0.
\]
It follows that $H(\blamb+n)+(-1)^rh_{n-r}\Delta_{(\lambda_1+1,\ldots.\lambda_r+1)}(H_r)=0$ for all $r+1\leq n\leq 2r-1$. For $n\geq 2r$ one 
uses
\[
H(\blamb+n)=e_1H(\blamb+n-1)-\ldots- (-1)^r e_rH(\blamb+n-r)
\] 
and induction,  to prove that $H(\blamb+n)-(-1)^{r}h_{n-r}\Delta_{(\lambda_1+1,\ldots.\lambda_r+1)}(H_r)=0$.
Therefore $G^\vee_r(z)\Delta_\blamb(H_r)$ involves no positive powers of $z$. Let us look now for the coefficients of negative powers. For  $0\leq n\leq \lambda_1$, the coefficient of $z^{-n}$ is the determinant:
\[
H(\blamb-n):=\left|\matrix{h_{\lambda_i+1-i-n}\cr h_{\lambda_i+2-i}\cr\vdots\cr h_{\lambda_i+r-i}}
\right|_{1\leq i\leq r}
\]
Therefore
\[
G_r^\vee(z)\Delta_\blamb(H_r)=\sum_{n=0}^{\lambda_1}{1\over z^n}H(\lambda-n)=\sum_{n=0}^{\lambda_1}\left|\matrix{\displaystyle{h_{\lambda_i+1-i-n}\over z^n}\cr h_{\lambda_i+2-i}\cr\vdots\cr h_{\lambda_i+r-i}}
\right|_{1\leq i\leq r}=\left|\matrix{\displaystyle{\sum_{n=0}^{\lambda_1}}\displaystyle{h_{\lambda_i+1-i-n}\over z^n}\cr h_{\lambda_i+2-i}\cr\vdots\cr h_{\lambda_i+r-i}}
\right|_{1\leq i\leq r }=
\]
\be
=\left|\matrix{G_r^\vee(z)h_{\lambda_1}&G_r^\vee(z)h_{\lambda_2-1}&\ldots&G_r^\vee(z)h_{\lambda_r-r+1}\cr
h_{\lambda_1+1}&h_{\lambda_2}&\ldots&h_{\lambda_r+r-2}\cr
\vdots&\vdots&\ddots&\vdots\cr
h_{\lambda_1+r-1}&h_{\lambda_2+r-2}&\ldots&h_{\lambda_r}}
\right|\label{eq:prematG}
\ee
Let $R_1,R_2,\ldots, R_r$ be the rows of the matrix~(\ref{eq:prematG}) and let
\[
G_r^\vee(z)R_j= (G_r^\vee(z)h_{\lambda_1+j-1},\ldots,G_r^\vee(z)h_{\lambda_j},\ldots,G_r^\vee(z)h_{\lambda_r+r-j}).
\]
Then  
\[
G_r^\vee(z)R_j=R_j+{1\over z}G_r^\vee(z)R_{j-1},
\]
 for all $2\leq j\leq r$. Therefore, again by the skew-symmetry and multi-linearity of the determinant, one obtains:
 \begin{eqnarray*}
 G_r^\vee(z)\Delta_\blamb(H_r)&=&\left|\matrix{G_r^\vee(z)h_{\lambda_1}&G_r^\vee(z)h_{\lambda_2-1}&\ldots&G_r^\vee(z)h_{\lambda_r-r+1}\cr
h_{\lambda_1+1}&h_{\lambda_2}&\ldots&h_{\lambda_r+r-2}\cr
\vdots&\vdots&\ddots&\vdots\cr
h_{\lambda_1+r-1}&h_{\lambda_2+r-2}&\ldots&h_{\lambda_r}}
\right|=\\ \\
&=& \left|\matrix{G_r^\vee(z)h_{\lambda_1}&G_r^\vee(z)h_{\lambda_2-1}&\ldots&G_r^\vee(z)h_{\lambda_r-r+1}\cr\cr
G_r^\vee(z)h_{\lambda_1+1}&G_r^\vee(z)h_{\lambda_2}&\ldots&G_r^\vee(z)h_{\lambda_r+r-2}\cr\cr\
\vdots&\vdots&\ddots&\vdots\cr\cr
G_r^\vee(z)h_{\lambda_1+r-1}&G_r^\vee(z)h_{\lambda_2+r-2}&\ldots&G_r^\vee(z)h_{\lambda_r}}
\right|=\Delta_\blamb(G_r^\vee(z)H_r)
 \end{eqnarray*}
and the Theorem is  proven.\qed
\bclm{\bf Corollary.} {\em For each $r\geq 0$, by abuse of notation, let $G_r(z), G^\vee_r(z):B_r[z^{-1}]\sra B_r[z^{-1}]$ be the $\QQ[z^{-1}]$ linear extension of the corresponding maps $B_r\sra B_r[z^{-1}]$. Then
\[
G_r(z)\circ G^\vee_r(z)=G^\vee_r(z)\circ G_r(z)=1_{B_r[z^{-1}]}
\]
i.e. they are inverse of each other.
}
\eclm
\proof
It suffices to evaluate each composition on $\Delta_\blamb(H_r)$:
\[
\hskip20ptG_r(z)\circ G^\vee_r(z)\Delta_\blamb(H_r)=G_r(z)\Delta_\blamb(G^\vee_r(z)H_r)=\Delta_\blamb(G_r(z)G^\vee_r(z)H_r)=\Delta_\blamb(G_r(z)).\hskip30pt \qed
\]
%\begin{eqnarray*}
%\hskip40pt G_r(z)\circ G^\vee_r(z)\Delta_\blamb(H_r)&=&G_r(z)\Delta_\blamb(G^\vee_r(z)H_r)=\\
%&=&\Delta_\blamb(G_r(z)G^\vee_r(z)H_r)=\Delta_\blamb(G_r(z)).\hskip110pt \qed
%\end{eqnarray*}
\claim{\bf Remark.} Analogously to Remark~\ref{rmkdetG}, if one defines
\[
\widetilde{G}^\vee_r(z)u_j=\sum_{n\geq 0} G_r(z)h_{n+j}{t^n\over n!}=\sum_{i\geq 0}{u_{j-i}\over z^i},
\]
then $(G^\vee_r(z)\Delta_\blamb(H_r))\cdot \Phi^r_0=\det(\widetilde{G}^\vee_r(z))\cdot \Phi^r_{0,\blamb}=$
\begin{eqnarray*}
%(G^\vee_r(z)\Delta_\blamb(H_r))\cdot \Phi^r_0&=&\det(\widetilde{G}^\vee_r(z))\cdot \Phi^r_{0,\blamb}=\\
&=&\widetilde{G}^\vee_r(z)u_{\lambda_1}\w \widetilde{G}^\vee_r(z) u_{-1+\lambda_2}\w\ldots\w\widetilde{G}^\vee_r (z)u_{-r+1+\lambda_r}\w \Phi^{r}_{-r}.
\end{eqnarray*}

\bclm{\bf Corollary.} 
\label{qucomv} {\em Let $h_{i_1}\cdot \ldots\cdot h_{i_s}$ be an arbitrary product of terms of $H_r$, with $s\leq r$. Then
\[
G^\vee_r(z)(h_{i_1}\cdot\ldots\cdot h_{i_s})=G_r(z) h_{i_1}\cdot\ldots\cdot G_r(z)h_{i_s}
\] 
}
\eclm
\proof Using the same argument, mutatis mutandis, as in Corollary~\ref{qucom}.\qed

%\bclm{\bf Corollary.} {\em For each $n\geq 0$
%\[
%G^\vee_r(z) h_n=\sum_{0\leq i\leq n}{h_{n-i}\over z^i}=h_n+{h_{n-1}\over z}+\ldots+{h_1\over z^{n-1}}+{1\over z^n}
%\]
%}
%\eclm

\bclm{\bf Corollary.}{ \em If $r=\infty$, then $G^\vee_\infty(z):B_r[z^{-1}]\sra B_r[z^{-1}]$ is a ring homomorphism.}\qed
\eclm

\section{The case $r=\infty$}\label{infty}
\claim{} We propose now an alternative proof for  the expressions of $\Gamma_\infty(z)$ and $\Gamma^\vee_\infty(z)$, with respect to that shown e.g. in~\cite{KR}, using the fact $G_\infty(z),G^\vee_\infty(z):B_\infty\sra B_\infty[z^{-1}]$ are ring homomorphisms. Recall that  $B_\infty=\QQ[H_\infty]=\QQ[X_\infty]$ and that the terms of the sequences $H_\infty=(h_1,h_2,\ldots)$ and $X_\infty=(x_1,x_2,\ldots)$ are algebraically independent in this case. Then relations~(\ref{eq:derhni}) hold.
Recall also that if $\Dcal(z)$ is a formal power series whose coefficients are first order differential operators on a $\QQ$-algebra $A$, then
\[
\exp(\Dcal(z)):A\sra A[[z]]
\]
is a  homomorphism of $\QQ$-algebras, in the sense that $\exp(\Dcal(z))(ab)=\exp(\Dcal(z))(a)\exp(\Dcal(z))(b)$.

\claim{} By Corollary~\ref{corhnz}:
\[
G_\infty(z)h_n=h_n-{h_{n-1}\over z}=\left(1-{1\over z}{\partial\over \partial x_1}\right)h_n
\]
where the last equality is because of~(\ref{eq:derhni}). Using the well known identity:
\[
1-{a}=\exp\left(-\sum_{n\geq 0} {a^n\over n}\right)
\]
for  $a=\displaystyle{1\over z}{\partial\over \partial x_1}$, one obtains:

\[
G_\infty(z)h_n=\exp \left(-\sum_{i\geq 1}{1\over iz^i}{\partial^n\over \partial x_1^n}\right)h_n=\exp \left(-\sum_{n\geq 1}{1\over iz^i}{\partial\over \partial x_i}\right)h_n.
\]
Since
\be
\Dcal\left({1\over z}\right)=\sum_{i\geq 1}{1\over iz^i}{\partial\over \partial x_i}\label{eq:d1ovz}
\ee
 is a formal power series in the indeterminate $1/z$ and the coefficients are the  first order differential operators $\partial/\partial x_i$, it follows that
$\exp (-\Dcal(1/z)):B_\infty\sra B_\infty[z^{-1}]$ is a ring homomorphism such that $G_\infty(z)h_n=\exp (-D(1/z))h_n$.  Then $G_\infty(z)=\exp(-\Dcal(1/z))$, because $(h_i)_{i\geq 1}$ generate $B_\infty$ as  a $\QQ$-algebra. In conclusion:
\[
\Gamma(z)={1\over E_\infty(z)}G_\infty(z)=\exp\left(\sum_{i\geq 1}x_it^i\right)\cdot \exp \left(-\sum_{n\geq 1}{1\over nz^n}{\partial\over \partial x_n}\right)
\]
which is precisely  expression~\cite[formula 5.25a]{KR} for $m=-1$, up to a factor used to keep track that $X(z)\w$, which is defined on $F^\infty_i$ for all $i\in\ZZ$,  is currently  operating on  $F_{-1}^\infty$.

The same kind of argument works for $G_\infty^\vee(z)$. One uses now the identity
\be
{1\over 1-\displaystyle{t\over z}}=1+\sum_{n\geq 1}{t^n\over z^n}=\exp\left(\sum_{n\geq 1}{t^n\over nz^n}\right).\label{eq:hidizi}
\ee
Thus
\begin{eqnarray*}
G_\infty^\vee(z)h_n=\sum_{i\geq 0}{h_{n-i}\over z^i}=\left(1+\sum_{i\geq 1}{1\over z^i}{\partial^i\over \partial x_1^i}\right)h_n
\end{eqnarray*}
Using~(\ref{eq:hidizi}):
\[
G_\infty^\vee(z)h_n=\exp\left(\sum_{i\geq 1}{1\over iz^i}{\partial^i\over\partial x_1^i}\right)h_n=\exp\left(\sum_{i\geq 1}{1\over iz^i}{\partial\over\partial x_i}\right)h_n.
\]
again by virtue of~(\ref{eq:derhni}).
Hence, for each $n\geq 0$, $G_\infty^\vee(z)h_n=\exp(\Dcal(1/z))h_n$, where $\exp(\Dcal(1/z))$ is as in~(\ref{eq:d1ovz}) . Then $G_\infty^\vee(z)=\exp(\Dcal(1/z))$,  because they are both algebra homomorphisms. In conclusion:
\[
\Gamma^\vee_\infty(z)={E_\infty\over z}\exp\left(\sum_{n\geq 1}{1\over nz^n}{\partial^n\over \partial x_n}\right)={1\over z}\exp(-\sum_{i\geq 1}x_iz^i)\exp\left(\sum_{n\geq 1}{1\over nz^n}{\partial\over \partial x_n}\right),
\]
which is precisely  expression~\cite[formula 5.25b]{KR} for $m=1$, up to a factor used to keep track that $X^\vee(z)\lrcorner$, which is defined on $F^\infty_i$, for all $i\in \ZZ$,   is currently operating on  $F_{1}^\infty$.

\bigskip
E-mail address: {\bf {\tt letterio.gatto@polito.it, parham@ibilce.unesp.br}}

\medskip
\sc{Dipartimento di Scienze Matematiche, Politecnico di Torino, Italy}

\smallskip
\sc{IBILCE-UNESP, Campus De S\~ao Jos\'e  do Rio Preto, Brazil}

\end{document}